\documentclass[11pt]{article} 
\usepackage[utf8]{inputenc}
\usepackage[T1]{fontenc}
\usepackage[english]{babel}

\usepackage{amsmath,amsfonts,amssymb,amsthm}
\usepackage{enumitem}
\usepackage{mathtools}
\usepackage{booktabs}

\usepackage[mono=false]{libertine}
\usepackage[cmintegrals,libertine]{newtxmath}
\usepackage[cal=euler, scr=boondoxo]{mathalfa}

\useosf
\usepackage[a4paper,vmargin={3.5cm,3.5cm},hmargin={2.5cm,2.5cm}]{geometry}
\usepackage[font={small,sf}, labelfont={sf,bf}, margin=1cm]{caption}
\usepackage[pdftex,colorlinks=true]{hyperref}
\usepackage[pdftex]{color,graphicx}

\usepackage{stackrel}
\usepackage{soul} 

\usepackage{ifthen}
\newlength{\leftstackrelawd}
\newlength{\leftstackrelbwd}
\def\leftstackrel#1#2{\settowidth{\leftstackrelawd}%
{${{}^{#1}}$}\settowidth{\leftstackrelbwd}{$#2$}%
\addtolength{\leftstackrelawd}{-\leftstackrelbwd}%
\leavevmode\ifthenelse{\lengthtest{\leftstackrelawd>0pt}}%
{\kern-.5\leftstackrelawd}{}\mathrel{\mathop{#2}\limits^{#1}}}

\theoremstyle{plain}
\newtheorem{theorem}{Theorem}

\newtheorem{proposition}[theorem]{Proposition}
\newtheorem{lemma}[theorem]{Lemma}
\theoremstyle{definition}

\newcommand{\ind}{\mathbf{1}}

\newcommand{\R}{\mathbb{R}}

\newcommand{\Z}{\mathbb{Z}}
\newcommand{\Q}{\mathbb{Q}}
\newcommand{\rmd}{\mathrm{d}}

\newcommand{\map}{\mathfrak{m}}

\begin{document}
	
\title{\bf On polynomials counting essentially irreducible maps}
\author{\textsc{Timothy Budd}\footnote{Radboud University, Nijmegen, The Netherlands. Email: \href{mailto:t.budd@science.ru.nl}{t.budd@science.ru.nl}}}
\date{\today}
\maketitle
\begin{abstract}
	We consider maps on genus-$g$ surfaces with $n$ (labeled) faces of prescribed even degrees.
	It is known since work of Norbury that, if one disallows vertices of degree one, the enumeration of such maps is related to the counting of lattice points in the moduli space of genus-$g$ curves with $n$ labeled points and is given by a symmetric polynomial $N_{g,n}(\ell_1,\ldots,\ell_n)$ in the face degrees $2\ell_1, \ldots, 2\ell_n$. 
	We generalize this by restricting to genus-$g$ maps that are essentially $2b$-irreducible for $b\geq 0$, which loosely speaking means that they are not allowed to possess contractible cycles of length less than $2b$ and each such cycle of length $2b$ is required to bound a face of degree $2b$.
	The enumeration of such maps is shown to be again given by a symmetric polynomial $\hat{N}_{g,n}^{(b)}(\ell_1,\ldots,\ell_n)$ in the face degrees with a polynomial dependence on $b$.
	These polynomials satisfy (generalized) string and dilaton equations, which for $g\leq 1$ uniquely determine them.
	The proofs rely heavily on a substitution approach by Bouttier and Guitter and the enumeration of planar maps on genus-$g$ surfaces. 
\end{abstract}	

\section{Introduction}

The enumeration of maps on surfaces (or \emph{fatgraphs} or \emph{ribbon graphs}) has a long history in mathematics and physics. 
Since the seminal results of Tutte in the sixties (including \cite{Tutte1962a}) many families of maps have been enumerated using a great variety of tools, including loop equations and analytic combinatorics, character theory of the symmetric group, random matrix theory and tree bijections.
Among these are many families corresponding to maps with restrictions on the face degrees and on the girth, i.e.\ the length of the shortest cycle, or the essential girth, which is the girth of the universal cover in the case of higher-genus maps.
For instance, loopless triangulations are maps with all faces restricted to be of degree three and (essential) girth at least two.
Until recently the constraints on the minimal girth were always rather small (2, 3 or 4).
This changed with the discovery by Bernardi \& Fusy of a tree encoding of planar maps with arbitrarily large minimal girth \cite{Bernardi2012a,Bernardi2012} which relies on the existence of canonical orientations on the edges of such maps (see also \cite{Albenque2015}).
Planar maps with the slightly stronger constraint of \emph{$d$-irreducibility}, meaning that they are required to have girth at least $d$ and that the only cycles of length $d$ are contours of a face of degree $d$, were enumerated by Bouttier \& Guitter in \cite{Bouttier2014,Bouttier2014a}.
Relying on a slice decomposition they also obtain a tree encoding, that is similar to the case with girth constraints.

The unified treatment of maps with girth or irreducibility constraint presents the opportunity to study the dependence of the map enumeration not only on the number of faces and their degrees but also on the girth parameter.
This is the goal of the current work.
As we will see this dependence is especially simple if we restrict our attention to maps with even face degrees and no vertices of degree one.
In particular, the number of $2b$-irreducible planar maps with $n$ labeled faces of even degrees $2\ell_1,\ldots,2\ell_n$ will be shown to depend polynomially on $b$ and $\ell_1,\ldots,\ell_n$.
The same is true for higher-genus maps when they are required to be essentially $2b$-irreducible, meaning that their universal cover is $2b$-irreducible.

The appearance of polynomials in the enumeration of maps with control on the face degrees, but without girth constraints, was already observed by Norbury in \cite{Norbury2010} and is closely linked to the decorated moduli space of genus $g$ curves with $n$ labeled points.
In the limit of large face degrees one can think of these maps as approximating the ribbon graphs with real edge lengths that appear in Kontsevich's proof \cite{Kontsevich1992} of Witten's conjecture. 
In particular, the coefficients of the leading order monomials in the face degrees encode intersection numbers of the first Chern classes of certain tautological bundles over moduli space.
The polynomials derived in this work for essentially $2b$-irreducible maps should contain more information since one may examine the leading order monomials not only in the face degrees but also in the irreducibility parameter $b$.
As will be demonstrated in the follow-up work \cite{Budd2020}, the homogeneous parts of top degree of these polynomials, seen as function of the face degrees as well as $b$, precisely compute certain volumes of ribbon graphs with real edge lengths and real girth constraints. 
We will observe that such volumes are in turn closely (but still mysteriously) related to the Weil-Petersson volumes of hyperbolic surfaces.

\subsection{Main results}

A \emph{genus-$g$ map} is a (multi)graph that is properly embedded in a surface of genus $g$, viewed up to orientation-preserving homeomorphisms of the surface.
Here \emph{properly embedded} means that edges only meet at their endpoints and that the complement of the graph is a disjoint union of topological disks.
We denote the set of vertices, edges, and faces of a map $\map$ by $\mathcal{V}(\map)$, $\mathcal{E}(\map)$ and $\mathcal{F}(\map)$ respectively.
A map is \emph{rooted} if it is equipped with a distinguished oriented edge, the \emph{root edge}.
Given a set $\mathcal{M}$ of maps with labeled faces and $\vec{\mathcal{M}}$ the corresponding set of rooted maps, we enumerate $\mathcal{M}$ via the formula
\begin{equation*}
\| \mathcal{M} \| = \sum_{\map\in\mathcal{M}} \frac{1}{|\operatorname{Aut}(\map)|} = \sum_{\map\in \vec{\mathcal{M}}} \frac{1}{2|\mathcal{E}(\map)|},
\end{equation*}
where $\operatorname{Aut}(\map)$ is the group of orientation-preserving automorphisms of the map $\map$ preserving the face labels.
In many cases there are no non-trivial automorphisms, for instance when the maps are planar and have three or more faces.
For a family $\mathcal{M}$ of such maps we thus have that $\|\mathcal{M}\|$ is simply the cardinality of $\mathcal{M}$.

\begin{figure}
	\centering
	\includegraphics[width=.75\linewidth]{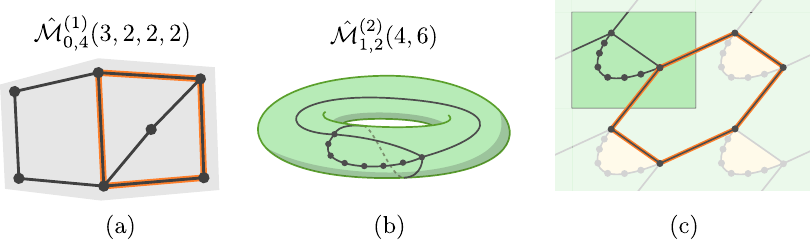}
	\caption{(a) Example of a planar map that is $2$-irreducible but not $4$-irreducible because it has a simple cycle of length $4$ (indicated in orange) that does not bound a face. (b) A genus-$1$ map that is essentially $4$-irreducible but not essentially $6$-irreducible, because its universal cover (c) has a simple cycle of length $6$ that does not bound a face. \label{fig:irreducibleexamples}}
\end{figure}

A planar map (or genus-$0$ map) is called \emph{$2b$-irreducible} \cite{Bouttier2014} for $2b\geq 0$ if it contains no simple cycle of length smaller than $2b$, i.e.\ it has \emph{girth} at least $2b$, and every simple cycle of length $2b$ is the boundary of a face of degree $2b$ (see Figure \ref{fig:irreducibleexamples}a).
By construction every map is $0$-irreducible.
A genus-$g$ map with $g\geq 1$ is said to be \emph{essentially $2b$-irreducible} if its universal cover viewed as an infinite planar map is $2b$-irreducible (Figure \ref{fig:irreducibleexamples}b \& c).

\begin{figure}
	\centering
	\includegraphics[width=.75\linewidth]{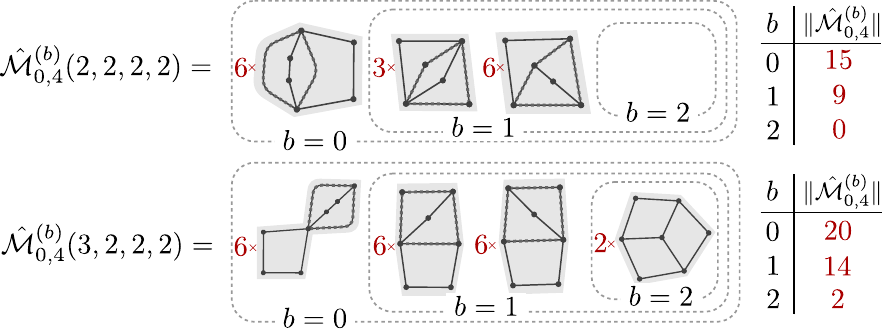}
	\caption{Illustration of the enumeration of $2b$-irreducible planar maps with four faces of degrees $2\ell_1,2\ell_2,2\ell_3,2\ell_4$ for $(\ell_1,\ell_2,\ell_3,\ell_4) = (2,2,2,2)$ and $(3,2,2,2)$ respectively. Some simple cycles that are relevant for the irreducibility constraints are represented as dotted curves. The number of inequivalent face labelings is indicated in red. As one may check the enumeration matches the formula $\|\hat{\mathcal{M}}^{(b)}_{0,4}\| = \ell_1^2+\ell_2^2+\ell_3^2+\ell_4^2 - (3b^2+3b+1) + 3\,\ind_{\{\ell_1=\cdots=\ell_4=b\}}$. \label{fig:mapsfourfaces}}
\end{figure}

For $g,b\geq 0$ and $n \geq 1$ (provided $n\geq 3$ if $g=0$) and $\ell_1, \ldots \ell_n \geq \max(b,1)$ we denote by $\mathcal{M}_{g,n}^{(b)}(\ell_1,\ldots,\ell_n)$ the set of essentially $2b$-irreducible genus-$g$ maps with $n$ labeled faces of degrees $2\ell_1,\ldots,2\ell_n$.
Let $\hat{\mathcal{M}}^{(b)}_{g,n}(\ell_1,\ldots,\ell_n)$ be the subset of such maps that contain no vertex of degree $1$.
Our first result is that $\|\hat{\mathcal{M}}^{(b)}_{g,n}(\ell_1,\ldots,\ell_n)\|$ is polynomial in the face degrees, apart from a small correction in the planar case when all degrees are equal to $2b$ (see Figure \ref{fig:mapsfourfaces} for an example).

\begin{theorem}\label{thm:mainpolynomial}
	For every $g,b\geq 0$ and $n\geq 1$ (provided $n\geq 3$ if $g=0$) there exists a symmetric polynomial $\hat{N}^{(b)}_{g,n}(\ell_1,\ldots,\ell_n)$ of degree $n+3g-3$ in $\ell_1^2, \ldots, \ell_n^2$ such that for all $\ell_1, \ldots \ell_n \geq \max(b,1)$,
	\begin{equation}\label{eq:mainenumeration}
	\|\hat{\mathcal{M}}^{(b)}_{g,n}(\ell_1,\ldots,\ell_n)\| = \hat{N}_{g,n}^{(b)}(\ell_1, \ldots, \ell_n) \,+\, \ind_{\{g=0,\,n\geq 4,\,\ell_1=\cdots=\ell_n=b\}}\,  \frac{(n-1)!}{2}(-1)^n
	\end{equation}
	and
	\begin{equation}\label{eq:MvsMhat}
	\|\mathcal{M}_{g,n}^{(b)}(\ell_1,\ldots,\ell_n)\| = \sum_{p_1=b}^{\ell_1}A^{(b)}_{\ell_1,p_1} \cdots \sum_{p_n=b}^{\ell_n}A^{(b)}_{\ell_n,p_n} \,\|\hat{\mathcal{M}}_{g,n}^{(b)}(p_1,\ldots,p_n)\|,
	\end{equation}
	where $A^{(b)}_{\ell,p}=\ind_{\ell=p=b} + \frac{p}{\ell} \binom{2\ell}{\ell-p}\,\ind_{\ell \geq p > b}$.
	Moreover, for fixed $g$ and $n$, the dependence of $\hat{N}^{(b)}_{g,n}(\ell_1,\ldots,\ell_n)$ on $b,\ell_1,\ldots,\ell_n$ is polynomial of degree $2n+6g-6$.
\end{theorem}

\noindent
Examples of the polynomials for small $g$ and $n$ are listed in Table~\ref{tab:polynomials}.
Next we show that the polynomials $\hat{N}_{g,n}^{(b)}$ and $\hat{N}_{g,n+1}^{(b)}$ satisfy two linear equations, which in the case of genus $0$ and $1$ completely characterize them in a recursive fashion.

\begin{theorem}\label{thm:mainstringdilaton}
	For every $g,b\geq 0$ and $n\geq 1$ (provided $n\geq 3$ if $g=0$) the polynomials of Theorem \ref{thm:mainpolynomial} satisfy the ``string equation''
	\begin{equation}\label{eq:string} 
	\hat{N}^{(b)}_{g,n+1}(\ell_1,\ldots,\ell_n,1) = \sum_{j=1}^n \sum_{k=b+1}^{\ell_j} 2 k\, \hat{N}^{(b)}_{g,n}(\ell_1,\ldots,\ell_{j-1},k,\ell_{j+1},\ldots,\ell_n) - \sum_{j=1}^n \ell_j \hat{N}^{(b)}_{g,n}(\ell_1,\ldots,\ell_n)
	\end{equation}
	and the ``dilaton equation''
	\begin{equation}\label{eq:dilaton} 
	\hat{N}_{g,n+1}^{(b)}(\ell_1,\ldots,\ell_n,1) - \hat{N}_{g,n+1}^{(b)}(\ell_1,\ldots,\ell_n,0) = (n+2g-2)\, \hat{N}_{g,n}^{(b)}(\ell_1,\ldots,\ell_n).
	\end{equation}
	In the planar and toroidal case these equations together with $\hat{N}_{0,3}^{(b)}(\ell_1,\ell_2,\ell_3)=1$ and $\hat{N}_{1,1}^{(b)}(\ell_1) = \frac{1}{12}\ell_1^2 - \frac{1}{12}$ uniquely determine the symmetric polynomials $\hat{N}_{g,n}^{(b)}$ for $g=0,1$ and all $n$.
	In general $\hat{N}^{(b)}_{g,1}(\ell_1)$ is independent of $b$.
\end{theorem}

\begin{table}
\begin{tabular*}{\textwidth}{@{}ccl@{}} \toprule
	$g$ & $n$ & $\hat{N}_{g,n}^{(b)}(\ell_1, \ldots, \ell_n)$ \\ \midrule
	$0$ & $3$ & $1$ \\
	& $4$ & $m_{(1)} - (3b^2 + 3b + 1)$ \\
	& $5$ & $\frac{1}{2} m_{(2)} + 2 m_{(1,1)}-\left(6 b^2+6 b+\frac{5}{2}\right) m_{(1)}+ (10 b^4+20 b^3+20 b^2+10 b+2)$ \\
	& $6$ & $\frac{1}{6}m_{(3)}+\frac{3}{2} m_{(2,1)}+6 m_{(1,1,1)}-\left(5 b^2+5 b+\frac{7}{3}\right) m_{(2)}-\left(18 b^2+18 b+9\right) m_{(1,1)}$\\
	& & \quad $+\left(30 b^4+60 b^3+65 b^2+35 b+\frac{49}{6}\right)
	m_{(1)}-\left(\frac{215 b^6}{6}+\frac{215 b^5}{2}+\frac{1085 b^4}{6}+\frac{365 b^3}{2}+\frac{340 b^2}{3}+40 b+6\right)$ \\ \midrule
	$1$ & $1$ & $\frac{1}{12}m_{(1)} - \frac{1}{12}$ \\
	& $2$ & $\frac{1}{24}m_{(2)}+\frac{1}{12} m_{(1,1)}-\frac{1}{8}m_{(1)}+\left(-\frac{b^4}{24}-\frac{b^3}{12}+\frac{b^2}{24}+\frac{b}{12}+\frac{1}{12}\right)$\\
	& $3$ & $\frac{1}{72}m_{(3)}+\frac{1}{12} m_{(2,1)}+\frac{1}{6} m_{(1,1,1)}-\left(\frac{b^2}{24}+\frac{b}{24}+\frac{1}{9}\right)
	m_{(2)}-\frac{1}{3} m_{(1,1)}+\left(-\frac{b^4}{12}-\frac{b^3}{6}+\frac{b^2}{8}+\frac{5 b}{24}+\frac{19}{72}\right) m_{(1)}$\\
	& & \quad $+\left(\frac{b^6}{18}+\frac{b^5}{6}+\frac{2
		b^4}{9}+\frac{b^3}{6}-\frac{5 b^2}{18}-\frac{b}{3}-\frac{1}{6}\right)$\\ \midrule
	$2$ & $1$ & $\frac{1}{6912}m_{(4)}-\frac{13}{5760}m_{(3)}+\frac{119 }{11520}m_{(2)}-\frac{143}{8640}m_{(1)}+\frac{1}{120}$ \\
	& $2$ & $\frac{1}{34560}m_{(5)}+\frac{29 }{17280}m_{(3,2)}+\frac{1}{2304}m_{(4,1)}-\frac{1}{1280}m_{(4)}-\frac{317
		}{17280}m_{(2,2)}-\frac{73}{8640}m_{(3,1)}$\\
	& & \quad $+\frac{17}{2304}m_{(3)}+\frac{61}{1280}m_{(2,1)}-\frac{1009}{34560}m_{(2)}-\frac{1543 }{17280}m_{(1,1)}+\frac{137}{2880}m_{(1)}$ \\
	& & \quad $+ \left(-\frac{b^{10}}{34560}-\frac{b^9}{6912}+\frac{b^8}{2880}+\frac{13 b^7}{5760}-\frac{7 b^6}{11520}-\frac{119 b^5}{11520}-\frac{31 b^4}{17280}+\frac{143
		b^3}{8640}+\frac{b^2}{480}-\frac{b}{120}-\frac{1}{40}\right)$\\
	 \bottomrule
\end{tabular*}
\caption{The first few polynomials $\hat{N}_{g,n}^{(b)}(\ell_1,\ldots,\ell_n)$. They are expressed in terms of the basis $m_{(\alpha_1,\ldots,\alpha_p)}(\ell_1,\ldots,\ell_n)=\sum_{(\beta_1,\ldots,\beta_n)} \ell_1^{2\beta_1}\ell_2^{2\beta_2}\cdots\ell_n^{2\beta_n}$ of symmetric even polynomials, where the sum runs over permutations $(\beta_1,\ldots,\beta_n)$ of $(\alpha_1,\ldots,\alpha_p,0,\ldots,0)$. For example, $m_{(1)}(\ell_1)=\ell_1^2$, $m_{(2)}(\ell_1,\ell_2) = \ell_1^4+\ell_2^4$, $m_{(1,1)}(\ell_1,\ell_2,\ell_3) = \ell_1^2\ell_2^2+\ell_1^2\ell_3^2+\ell_2^2\ell_3^2$ and $m_{(3,1)}(\ell_1,\ell_2,\ell_3) = \ell_1^6\ell_2^2+\ell_1^6\ell_3^2+\ell_2^6\ell_1^2+\ell_2^6\ell_3^2+\ell_3^6\ell_1^2+\ell_3^6\ell_2^2$. \label{tab:polynomials}}
\end{table}

These theorems generalize results of Norbury in \cite{Norbury2010,Norbury2013} that deal with the case $b=0$ albeit in a broader setting where the faces degrees are allowed to be odd as well.
Remarkably the string and dilaton equations for essentially $2b$-irreducible maps differ from their $b=0$ counterpart only in the lower bound of the inner sum in the string equation \eqref{eq:string}.
The way we arrive to these equations, however, is quite different from the methods used by Norbury, which are obtained with the help of a topological recursion formula satisfied by genus-$g$ maps without vertices of degree one.
No such recursion formula is known when an irreducibility constraint is present.
Instead, we rely heavily on the substitution approach by Bouttier \& Guitter \cite{Bouttier2014,Bouttier2014a}, which generalizes the approach of Tutte \cite{Tutte1962a} and Mullin \& Schellenberg \cite{Mullin1968} for triangulations respectively quadrangulations with an irreducibility constraint to arbitrary maps with controlled face degrees.
Using known results on the enumeration of genus-$g$ maps it allows us to obtain sufficiently manageable generating functions for essentially $2b$-irreducible maps.
A cautionary note: the results of Theorem \ref{thm:mainpolynomial} and \ref{thm:mainstringdilaton} in the case $g\geq 1$ depend on the enumeration of genus-$g$ maps in terms of certain moments, summarized in Section \ref{sec:highergenuspartition}, that stem from topological recursion \cite{Eynard2016} or matrix models \cite{Ambjorn1993,Akemann1996}.
We leave it to the reader to judge the level of rigor of these sources.

Restricting the face degrees to be strictly larger than $2b$ in Theorem \ref{thm:mainpolynomial} and \ref{thm:mainstringdilaton} leads to analogous results for genus-$g$ maps with a constraint on the essential girth, i.e.\ the girth of the universal cover. 
In particular, if $\hat{\mathcal{G}}_{g,n}^{(\geq b)}(\ell_1,\ldots,\ell_n)$ denotes the set of genus-$g$ maps with $n$ labeled faces ($n\geq 3$ if $g=0$) of degrees $2\ell_1,\ldots,2\ell_n \geq 2b$, no vertices of degree one, and essential girth at least $2b > 0$, then
\begin{equation*}
\|\hat{\mathcal{G}}_{g,n}^{(\geq b)}(\ell_1,\ldots,\ell_n)\| = \hat{N}_{g,n}^{(b-1)}(\ell_1,\ldots,\ell_n)
\end{equation*}
is polynomial in $b,\ell_1,\ldots,\ell_n$.
If $\hat{\mathcal{G}}_{g,n}^{(b)}(\ell_1,\ldots,\ell_n)$ is the corresponding set with essential girth exactly $2b$ then
\begin{equation*}
\|\hat{\mathcal{G}}_{g,n}^{(b)}(\ell_1,\ldots,\ell_n)\| = \hat{N}_{g,n}^{(b-1)}(\ell_1,\ldots,\ell_n) - \ind_{\{\ell_1,\ldots,\ell_n > b\}}\hat{N}_{g,n}^{(b)}(\ell_1,\ldots,\ell_n),
\end{equation*}
which is again polynomial when restricted to $\ell_1,\ldots,\ell_n > b$.

\subsection{Questions}

Our main results are easily stated and very much analogous to the enumeration of maps without irreducibility constraint, but the proofs rely on some tedious computations.
This makes one suspect that a simpler combinatorial understanding is yet to be discovered.
Here we list several natural open questions in this direction.

\begin{enumerate}[label={\bfseries \arabic*.}]
	\item {\bfseries Do Theorem \ref{thm:mainpolynomial} and Theorem \ref{thm:mainstringdilaton} generalize to essentially $2b$-irreducible maps including odd face degrees?} 
	The methods we used to prove the results based on generating function should in principle work when odd face degrees are included, but the computations would become significantly more complicated.
	In the case $b=0$ it is known from the work of Norbury \cite[Theorem 1]{Norbury2010} that general maps without vertices of degree $1$ are enumerated by polynomials in the boundary length, where the polynomials depends on $g$ and $n$ as well as on the parity of the face degrees (i.e.\ the number of odd face degrees).
	\item {\bfseries Does the enumeration of essentially $2b$-irreducible maps admit a topological recursion?}\\
	It is far from clear whether essentially $2b$-irreducible maps admit a combinatorial decomposition like the one underlying Tutte's equation for general maps, because most natural decompositions one may consider do not preserve the irreducibility constraint.
	As alluded to above, in the case $b=0$ and general face degrees (even and odd) a recursion formula for the polynomials $N_{g,n}^{(0)}$ is known \cite[Theorem 4]{Norbury2010} and can be traced back to Tutte's equation.
	It is connected to the fact \cite[Theorem 2]{Norbury2013} that the polynomials appear as the coefficients of the Eynard-Orantin invariants of the plane curve $x y-y^2 = 1$.
	The string and dilaton equations are naturally associated to this curve \cite{Eynard2007,Norbury2013}.
	Given the resemblance to the string and dilaton equations for essentially $2b$-irreducible maps, one might hope to guess a plane curve whose invariants encode the polynomials $\hat{N}^{(b)}_{g,n}$.
	In \cite{Norbury2013a} the equations associated to a fairly broad family of plane curves were determined, but unfortunately they do not include the specific string and dilaton equations from Theorem \ref{thm:mainstringdilaton}.
	\item {\bfseries Do the string and dilaton equation have a combinatorial explanation?} In the case $b=0$ one can give a combinatorial interpretation to  $\hat{N}_{g,n+1}^{(0)}(\ell_1,\ldots,\ell_n,1)$ and $\hat{N}_{g,n+1}^{(0)}(\ell_1,\ldots,\ell_n,0)$ as the counting of certain maps with a distinguished face of degree $2$ or with a distinguished vertex \cite{Norbury2013}. 
	Such an interpretation is troublesome when $b > 1$, since Theorem \ref{thm:mainpolynomial} only holds for maps with face degrees that are $2b$ or greater. 
	Does there exist a generalized notion of essential $2b$-irreducibility that extends the validity of Theorem \ref{thm:mainpolynomial} to the case where one or several faces have degree less than $2b$? 
	\item {\bfseries Are there bijective interpretations to the generating functions of essentially $2b$-irreducible maps of genus $g\geq 1$?}
	As mentioned before planar maps with control on the girth can be studied via the existence of certain canonical orientations which in turn give rise to bijections with certain decorated trees \cite{Bernardi2012a,Bernardi2012,Albenque2015}.
	This has recently been extended to the toroidal case \cite{Fusy2020} (see also \cite{Bonichon2019} for a special case), where the encoding is via certain decorated unicellular toroidal maps.
	In the case of $2b$-irreducible planar maps it was shown in \cite[Section~6]{Bouttier2014} that an iterative decomposition into slices leads to  similar encodings by decorated trees.
	The higher-genus cases however are much less developed, due to a lack of good canonical orientations or slice decompositions.
\end{enumerate}

\subsection{Outline}

We start in Section \ref{sec:skeleton} by formulating a sufficient criterion for essential $2b$-irreducibility, that can easily be checked in a skeleton decomposition of a map. 
This turns the enumeration of essentially $2b$-irreducible maps into a problem of counting lattice points in convex polyhedra that are associated to the possible skeletons.
Although performing this counting in general is hard, it gives a glimpse of the appearance of polynomials in the enumeration, due to general results on the counting of lattice points in convex polyhedra.
In Section \ref{sec:genfun} we start afresh by summarizing the substitution approach of Bouttier \& Guitter and generalize it to arbitrary genus. 
With the help of enumeration results of arbitrary genus-$g$ maps with control on the (even) face degrees, we derive relatively succinct expressions for the generating functions of essentially $2b$-irreducible genus-$g$ maps (with vertices of arbitrary degree).
In Section \ref{sec:hatgenfun} we perform a further substitution to disallow vertices of degree one and derive Theorem \ref{thm:mainpolynomial} from general properties of the generating functions.
Finally, the string and dilaton equations are verified on the level of the generating functions in Section \ref{sec:stringdilaton}.

\subsection*{Acknowledgments}
This work is part of the START-UP 2018 programme with project number 740.018.017, which is financed by the Dutch Research Council (NWO). 
We warmly thank an anonymous referee for useful suggestions.

\section{A first glimpse of the enumeration}\label{sec:skeleton}

Let $\map$ be a map on an oriented closed surface $S$ of genus $g\geq 1$ and let $\pi : S^\infty \to S$ be the universal cover of $S$, such that $S^\infty$ has the topology of the plane.
Then there exists a unique infinite planar map $\map^\infty$ on $S^\infty$ corresponding to the lift of $\map$ along $\pi$.
With a slight abuse of notation we will also denote by $\pi$ the mappings $\pi : \mathcal{F}(\map^\infty) \to \mathcal{F}(\map)$ and $\pi : \mathcal{E}(\map^\infty) \to \mathcal{E}(\map)$ that send faces and edges of $\map^\infty$ to their counterpart in $\map$. 
A mapping $D : S^\infty \to S^\infty$ is called a \emph{deck transformation} if $\pi \circ D = \pi$. 
The set of all deck transformations forms a group under composition that is isomorphic to the fundamental group of $S$.
Again we will abuse notation by using the same symbol for the mappings $D : \mathcal{F}(\map^\infty) \to \mathcal{F}(\map^\infty)$ and $D : \mathcal{E}(\map^\infty) \to \mathcal{E}(\map^\infty)$ describing the permutation of faces and edges under the deck transformation.

Recall that we call $\map$ essentially $2b$-irreducible if $\map^\infty$ is $2b$-irreducible, following e.g.\ \cite[Section 2]{Bonichon2019}. 
This definition can be impractical, since it involves a criterion on the lengths of the infinitely many simple cycles of $\map^\infty$.
We start by formulating an equivalent criterion that only involves checking the lengths of finitely many paths in $\map$, the essentially simple cycles in $\map$.
We call a closed path in $\map$ an \emph{essentially simple cycle} if it lifts to a simple cycle in $\map^\infty$ that encloses at most one face in $\pi^{-1}(f) \subset\mathcal{F}(\map^\infty)$ for each face $f\in\mathcal{F}(\map)$.
Informally, it is a (non-backtracking, contractible) cycle on $\map$ that bounds a simply-connected region.

\begin{lemma}\label{lem:irrcriterion}
	For $g,b\geq 1$, a genus-$g$ map $\map$ with faces of degree at least $2b$ is essentially $2b$-irreducible if and only if each essentially simple cycle of $\map$ has length at least $2b$ with equality only if it is the contour of a face of degree $2b$.
\end{lemma}
\begin{proof}
	The latter condition is clearly necessary for the map to be $2b$-irreducible.
	In order to prove that it is sufficient, let $\map$ be a genus-$g$ map and let $2k\geq 2$ be its essential girth, i.e.\ the length of the shortest cycle in the universal cover $\map^\infty$ of $\map$.
	We claim that the simple cycles of length $2k$ in $\map^\infty$ are precisely the lifts of essentially simple cycles of length $2k$ in $\map$.
	Let us first see how this proves the lemma.
	Suppose that each essentially simple cycle of $\map$ has length at least $2b$ with equality only if it bounds a face of degree $2b$.
	Then we must have $k\geq b$. 
	If $k>b$, then $\map^\infty$ is clearly $2b$-irreducible because it has no simple cycles of length $2b$ or shorter.
	If $k=b$, the claim implies that each simple cycle of length $2b$ in $\map^\infty$ corresponds to the contour of a face, and therefore $\map^\infty$ is $2b$-irreducible.
	So in both cases $\map$ is essentially $2b$-irreducible.
	
	To establish the claim we follow \cite[Section 3.1]{Bouttier2014} and introduce $\mathcal{C}_k({\map}^\infty)$ to be the set of \emph{outermost} cycles of length $2k$ in ${\map}^\infty$, i.e.\ the cycles whose interior is not fully contained in a different cycle of length $2k$.	
	Up to deck transformations there are only finitely many cycles of length $2k$ on  ${\map}^\infty$, so that in particular there is a bound on the number of faces enclosed by any cycle of length $2k$ in ${\map}^\infty$.
	Since every cycle of length $2k$ is either outermost or its interior is contained in the strictly larger interior of another cycle of length $2k$, it follows that the interior of each cycle of length $2k$ in ${\map}^\infty$ is contained in that of at least one outermost cycle in $\mathcal{C}_k({\map}^\infty)$.
	In particular, $\mathcal{C}_k({\map}^\infty)$ is non-empty.
	
	In \cite[Section 3.1]{Bouttier2014} it was shown that two distinct outermost cycles in a planar map, like ${\map}^\infty$, cannot overlap, i.e. must have disjoint interiors.
	In particular, the interior of an outermost cycle $C\in\mathcal{C}_k({\map}^\infty)$ cannot contain two distinct lifts $f_1,f_2\in (\pi)^{-1}(f)$ of a single face $f$ of $\map$.
	To see this, let $D : S^\infty \to S^\infty$ be a deck transformation such that $D(f_1) = f_2$.
	Then $D(C)$ is an outermost cycle different from $C$ (if it were identical the deck transformation $D$ would have a finite orbit, but this cannot happen for a closed surface).
	If $C$ encircles $f_1$, then $D(C)$ encircles $f_2$, implying that $C$ cannot encircle $f_2$.
		
	It follows that each cycle in $\mathcal{C}_k({\map}^\infty)$ corresponds to a lift of an essentially simple cycle in $\map$.
	The same is then true for any simple cycle of length $2k$ in ${\map}^\infty$ that is not outermost, since its interior is contained in that of an outermost cycle. 
	This establishes the claim and concludes the proof.
\end{proof}

Applied to a genus-$g$ map $\map$ without vertices of degree one and all faces of degree at least $2b$, Lemma~\ref{lem:irrcriterion} states in the case $g\geq 1$ that $\map$ is essentially $2b$-irreducible if and only if each essentially simple cycle enclosing at least two faces has length larger than $2b$.
The requirements on the lengths of essentially simple cycles that enclose a single face are automatically satisfied.
If $\map$ is planar it is $2b$-irreducible if and only if each simple cycle that encloses at least two faces on both sides has length larger than $2b$.
In particular, the requirement to be essentially $2b$-irreducible is non-trivial only when the number of faces is $n\geq 2$ for $g\geq 1$ and $n\geq 4$ for $g=0$ (assuming face degrees to be at least $2b$).

\begin{figure}
	\centering
	\includegraphics[width=.85\linewidth]{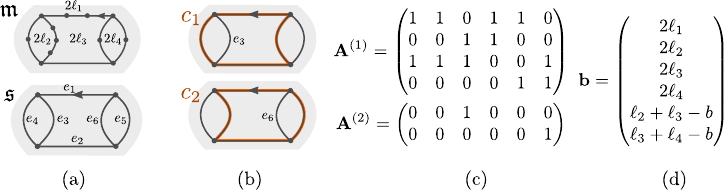}
	\caption{(a) A rooted genus-$0$ map $\map$ with four labeled faces and its skeleton $\mathfrak{s}=\mathsf{Skel}(\map)$ together with a labeling of its edges. (b) A complete collection $c_1,c_2$ of relevant simple cycles of $\mathfrak{s}$. (c) The face-edge incidence matrix $\mathbf{A}^{(1)}$ and the cycle-edge enclosure matrix $\mathbf{A}^{(2)}$ associated to $\mathfrak{s}$. (d) The vector $\mathbf{b}$ appearing on the right-hand side of the linear equation. \label{fig:skeleton}}
\end{figure}

A natural way to approach the enumeration of these maps is via the skeleton decomposition.
The \emph{skeleton} of a genus-$g$ map $\map$ without vertices of degree one is the genus-$g$ map $\mathsf{Skel}(\map)$ obtained by deleting each vertex of degree two and merging its incident edges (see Figure \ref{fig:skeleton}a).
In case $\map$ is rooted we take $\mathsf{Skel}(\map)$ to be rooted on the edge into which the root edge of $\map$ was merged (with the same orientation).
If we denote by $\hat{\vec{\mathcal{M}}}^{(b)}_{g,n}(\ell_1,\ldots,\ell_n)$ the set of rooted maps corresponding to $\hat{\mathcal{M}}^{(b)}_{g,n}(\ell_1,\ldots,\ell_n)$, we may express the enumeration of the latter as
\begin{align}
\|\hat{\mathcal{M}}^{(b)}_{g,n}(\ell_1,\ldots,\ell_n)\| &= \sum_{\map\in\hat{\vec{\mathcal{M}}}^{(b)}_{g,n}(\ell_1,\ldots,\ell_n)} \frac{1}{2|\mathcal{E}(\map)|}= \sum_{\map\in\hat{\vec{\mathcal{M}}}^{(b)}_{g,n}(\ell_1,\ldots,\ell_n)} \frac{\ind_{\{\text{root degree of }\map\text{ is at least }3\}} }{2|\mathcal{E}(\mathsf{Skel}(\map))|}\nonumber\\
&= \sum_{\mathfrak{s}} \frac{C^{(b)}_{\ell_1,\ldots,\ell_n}(\mathfrak{s})}{2|\mathcal{E}(\mathfrak{s})|},\label{eq:mhatskeleton}
\end{align}
where the last sum is over all rooted genus-$g$ maps $\mathfrak{s}$ with $n$ labeled faces and all vertices of degree at least $3$, 
\begin{equation}\label{eq:skelcount}
C^{(b)}_{\ell_1,\ldots,\ell_n}(\mathfrak{s}) = |\{ \map \in \hat{\vec{\mathcal{M}}}^{(b)}_{g,n}(\ell_1,\ldots,\ell_n) : \text{root degree at least }3\text{ and }\mathsf{Skel}(\map) = \mathfrak{s}\}|
\end{equation}
and the \emph{root degree} of $\map$ is the degree of the vertex at the origin of the root edge.
As we will demonstrate in the proposition below, $C^{(b)}_{\ell_1,\ldots,\ell_n}(\mathfrak{s})$ counts integer points in a certain convex polytope.

To understand this, we note that for fixed skeleton $\mathfrak{s}$ with edges $e_1,\ldots,e_{k}$, $k=|\mathcal{E}(\mathfrak{s})|$, the set of all rooted genus-$g$ maps $\map$ with root degree at least $3$ and $\mathsf{Skel}(\map) = \mathfrak{s}$ is in bijection with positive integer vectors $\mathbf{x}=(x_1,\ldots,x_k)\in \Z_{>0}^{k}$ by setting $x_i$ to be the number of edges of $\map$ that are merged into the edge $e_i$ of $\mathfrak{s}$.
Imposing that the $n$ faces of $\map$ have degrees $2\ell_1,\ldots,2\ell_n$ describes a set of $n$ linear equalities on $\mathbf{x}$.

It is not hard to see that the essential irreducibility constraint turns into a finite system of linear inequalities on $\mathbf{x}$.
Indeed, there is a one-to-one correspondence between the essentially simple cycles of $\map$ and those of its skeleton $\mathfrak{s}$.
Letting $c_1, \ldots, c_p$ be a complete list of essentially simple cycles of $\mathfrak{s}$ that enclose at least two faces (see Figure \ref{fig:skeleton}b), the map $\map$ is essentially $2b$-irreducible if and only if the integers $\mathbf{x}$ summed over the edges of each cycle $c_i$ exceed $2b$.
One may conveniently turn these inequalities on $\mathbf{x}$ into linear equalities by introducing auxiliary variables $x_{k+1},\ldots,x_{k+p}$, one for each cycle $c_p$.

We may then rely on quite general results (see e.g.\ \cite{DeLoera2005} for a gentle introduction) about counting integer points in polyhedra to learn something about $C^{(b)}_{\ell_1,\ldots,\ell_n}(\mathfrak{s})$.
In order to state the result in our case, we need to introduce some terminology. 
A function $\Z^k \to \Q$ is \emph{quasi-polynomial of degree $m$} if it can be expressed as a ($k$-variate) polynomial of degree $m$ with coefficients that are periodic functions $\Z^k\to \Q$.
A function $\Z^k \to \Q$ is \emph{piecewise quasi-polynomial of degree $m$} if one may subdivide $\Z^k$ into finitely many polyhedral regions such that restricted to each region it is quasi-polynomial and the maximal degree over all regions is $m$.

\begin{proposition}\label{thm:quasipolynomial}
$\|\hat{\mathcal{M}}^{(b)}_{g,n}(\ell_1,\ldots,\ell_n)\|$ is piecewise quasi-polynomial in $b,\ell_1,\ldots,\ell_n$ of degree $2n+6g-6$.
It is also of degree $2n+6g-6$ in $\ell_1,\ldots,\ell_n$ alone.
\end{proposition}
\begin{proof}
Let $\mathfrak{s}$ be a fixed genus-$g$ rooted map with $n$ labeled faces, $k$ edges and no vertices of degree less than three, and $c_1, \ldots,c_p$ a complete list of essentially simple cycles as above.
In the case $g=0$ we assume one side of each cycle $c_i$ is designated as the interior (the result will not depend on this choice).
When $g\geq 1$ we call the interior of $c_i$ the simply-connected region surrounded by $c_i$.
We claim that $C^{(b)}_{\ell_1,\ldots,\ell_n}(\mathfrak{s})$ takes the form of a \emph{vector partition function}, meaning that it counts positive integer solutions to a linear equation,
\begin{equation}\label{eq:vectorpartition}
C^{(b)}_{\ell_1,\ldots,\ell_n}(\mathfrak{s}) = |\{ \mathbf{x} \in \Z_{>0}^{k+p} : \mathbf{A} \mathbf{x} = \mathbf{b} \}|,
\end{equation}
where $\mathbf{A}$ is an $(n+p)\times(k+p)$ matrix of full rank with non-negative integer entries depending only on $\mathfrak{s}$, while $\mathbf{b} \in \Z_{>0}^{n+p}$ depends linearly on $b,\ell_1,\ldots,\ell_n$.
This puts us precisely in the setting of \cite[Theorem 1]{Sturmfels1995}\footnote{Note that we require $\mathbf{x}>0$ while the theorem applies to the situation $\mathbf{x}\geq 0$, but this difference amounts to a shift by one in $\mathbf{b}$.} (see also \cite[Theorem 2]{Beck2004}) which states that $C^{(b)}_{\ell_1,\ldots,\ell_n}(\mathfrak{s})$ is a piecewise quasi-polynomial in $b,\ell_1,\ldots,\ell_n$ of degree $k-n$.
By \eqref{eq:mhatskeleton}, $\|\hat{\mathcal{M}}^{(b)}_{g,n}(\ell_1,\ldots,\ell_n)\|$ amounts to a finite sum of such functions $C^{(b)}_{\ell_1,\ldots,\ell_n}(\mathfrak{s})$, so it is piecewise quasi-polynomial as well.
Its degree is equal to $k-n$ where $k$ is the maximal number of edges in a skeleton.
This maximum is achieved when all vertices are of degree three, in which case $k=3n+6g-g$.
This gives the claimed statement.

It remains to prove \eqref{eq:vectorpartition}.
Let $\mathbf{A}^{(1)}$ be the face-edge incidence matrix of dimension $n \times k$, meaning that $\mathbf{A}^{(1)}_{ij}$ is determined by whether the edge $e_j$ is incident on both sides to the $i$th face ($\mathbf{A}^{(1)}_{ij}=2$), on one side ($\mathbf{A}^{(1)}_{ij}=1$), or is not incident ($\mathbf{A}^{(1)}_{ij}=0$).
The matrix $\mathbf{A}^{(2)}$ of dimension $p \times k$ indicates which edges are surrounded by the cycles $c_i$. 
More precisely, we set $\mathbf{A}^{(2)}_{ij}=1$ if the edge $e_j$ is not part of the cycle $c_i$ but is contained in its interior and $\mathbf{A}^{(2)}_{ij}=0$ otherwise.
Let ${\boldsymbol \ell} \in \Z^n$ be the vector of face half-degrees $\ell_1, \ldots, \ell_n$. 
Finally the vector $\mathbf{m}=(m_1,\ldots,m_p)\in\Z^p$ is defined by setting $m_i = -b + \sum \ell_j$ where the sum is over the faces in the interior of the cycle $c_i$.
See Figure \ref{fig:skeleton}c for an example.
The matrix $\mathbf{A}$ and the vector $\mathbf{b}$ are given by
\begin{equation*}
 \quad \mathbf{A} = \begin{pmatrix}
\mathbf{A}^{(1)} & 0 \\
\mathbf{A}^{(2)} & \mathbb{I}_p
\end{pmatrix}, \quad \mathbf{b} = \begin{pmatrix}
2{\boldsymbol \ell}\\
\mathbf{m}
\end{pmatrix}.
\end{equation*}
Under the bijection described above between $(x_1,\ldots,x_k)\in \Z_{>0}^k$ and maps $\map$ with skeleton $\mathfrak{s}$, it is clear that $\sum_{j=1}^k\mathbf{A}^{(1)}_{ij}x_j = 2 \ell_i$ implements the face constraints correctly.
In terms of $(x_1,\ldots,x_k)$, the length $L_i$ of the cycle in $\map$ corresponding to $c_i$ is given by $L_i=\sum \mathbf{x}_j$ where the sum is over the edges $e_j$ in $c_i$.
Computing the sum of the degrees of the faces in the interior of $c_i$ in two ways leads to the identity $L_i + 2\sum_{j=1}^E\mathbf{A}^{(2)}_{ij}\mathbf{x}_j= 2b + 2\mathbf{m}_i$. 
The map $\map$ is essentially $2b$-irreducible precisely if $L_i > 2b$ for each $i=1,\ldots,p$.
This is clearly equivalent to the existence of $x_{k+1},\ldots,x_{k+p}\in \Z_{>0}$ such that $x_{k+i} + \sum_{j=1}^k\mathbf{A}^{(2)}_{ij}x_j = \mathbf{m}_i$.
Hence \eqref{eq:vectorpartition} really counts the maps in \eqref{eq:skelcount}.
\end{proof}

In principle, we could prove Theorem \ref{thm:mainpolynomial} for any fixed choice of $g$ and $n$ by tabulating the finitely many genus-$g$ skeleton maps with $n$ faces, solving the vector partition function \eqref{eq:vectorpartition} for each skeleton and observe that miraculously the piecewise quasi-polynomials add up to a regular polynomial.
As an illustration, let us look at $g=0$ and $n=3$, in which case the irreducibility does not play a role. 
There are precisely $7$ different labeled skeleton maps with three labeled faces (modulo rerooting) and for each of them $C^{(b)}_{\ell_1,\ell_2,\ell_3}(\mathfrak{s})$ is an indicator function.
Adding them all up we indeed find
\begin{align*}
\|\hat{\mathcal{M}}^{(b)}_{0,3}(\ell_1,\ell_2,\ell_3)\| &= \ind_{\{\ell_3 < \ell_1+\ell_2,\,\ell_2 < \ell_1+\ell_3,\,\ell_1 < \ell_2+\ell_3\}}+ \ind_{\{\ell_3 > \ell_1+\ell_2\}} + \ind_{\{\ell_2 > \ell_1+\ell_3\}} + \ind_{\{\ell_1 > \ell_2+\ell_3\}} \\ 
&\quad + \ind_{\{\ell_3 = \ell_1+\ell_2\}} + \ind_{\{\ell_2 = \ell_1+\ell_3\}} + \ind_{\{\ell_1 = \ell_2+\ell_3\}} = 1.
\end{align*}
Even when implemented on a computer the same exercise becomes intractable already for moderately large $n$ or $g$.

\section{Generating functions of essentially $2b$-irreducible maps}\label{sec:genfun}

For $b\geq 0$ and $g\geq 0$ let $\mathcal{M}^{(b)}_{g}$ be the set of essentially $2b$-irreducible genus-$g$ maps with faces of arbitrary even degree and $\vec{\mathcal{M}}^{(b)}_{g}$ the set of such maps that are rooted.
For our purpose we exclude planar maps with one or two faces from $\mathcal{M}^{(b)}_0$ and $\vec{\mathcal{M}}^{(b)}_0$.
For $b\geq 1$, the $2b$-irreducible genus-$g$ \emph{partition function} is given by the formal generating function
\begin{equation}\label{eq:partitionfunction}
F_g^{(b)}(x_{b},x_{b+1}, \ldots) = \sum_{\map\in\mathcal{M}^{(b)}_g} \frac{1}{|\operatorname{Aut}(\map)|} \prod_{f\in \mathcal{F}(\map)} x_{\deg(f)/2}  = \sum_{\map\in\vec{\mathcal{M}}^{(b)}_g} \frac{1}{2|\mathcal{E}(\map)|} \prod_{f\in \mathcal{F}(\map)} x_{\deg(f)/2}.
\end{equation}
For $b=0$ we set
\begin{equation}
F_g^{(0)}(x_0,x_1,\ldots) = \sum_{\map\in\mathcal{M}^{(0)}_g} \frac{1}{|\operatorname{Aut}(\map)|} (1+x_0)^{|\mathcal{V}(\map)|} \prod_{f\in \mathcal{F}(\map)} x_{\deg(f)/2}
\end{equation}
and we let $F_g(x_1,x_2\ldots) = F_g^{(0)}(0,x_1,x_2,\ldots)$ be the standard generating function of maps without vertex weights.

Even though we are working with an infinite number of formal generating variables $x_b,x_{b+1},\ldots$, there is no need to pay attention to convergence issues.
In fact, we could set an upper limit $2d$ on the face degrees and only consider generating functions in the variables $x_b,\ldots,x_d$ and it would not affect the expressions and proofs as long as $d>b$.
Since this upper limit will not play any role in the results while introducing clutter in the exposition, we choose to omit it in the following.
In any case we have (for $b \leq \ell_i \leq d$) that
\begin{equation}\label{eq:numfrompartition}
\|\mathcal{M}^{(b)}_{g,n}(\ell_1,\ldots,\ell_n)\| = \frac{\partial^n F_g^{(b)}}{\partial x_{\ell_1}\cdots\partial x_{\ell_n}}(0,0,\ldots).
\end{equation}

The goal of this section is to obtain manageable expressions for the partition functions.
We start with the genus-$0$ case, where the enumeration problem has been largely solved by Bouttier and Guitter in \cite{Bouttier2014,Bouttier2014a} using a \emph{substitution approach} that we summarize here.

\subsection{Substitution approach in planar case} 

\begin{figure}
	\centering
	\includegraphics[width=.8\linewidth]{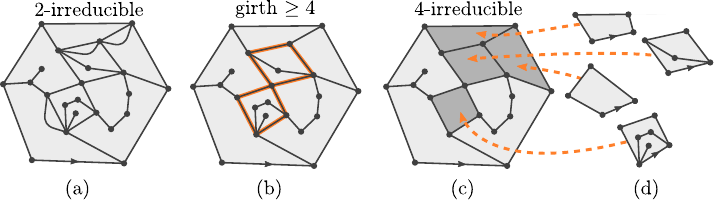}
	\caption{Example in the case $b=2$ and $\ell=3$. (a) A $(2b-2)$-irreducible rooted map in $\vec{\mathcal{M}}_{0,\ell}^{(b-1)}$. (b) A rooted map of girth at least $2b$ in $\mathcal{G}_{0,\ell}^{(\geq b)}$ is a $(2b-2)$-irreducible map without faces of degree $2b-2$. The outermost cycles of length $2b$ (excluding the faces of degree $2b$) are highlighted in orange. (c) Removing the interior of the outermost cycles gives a $2b$-irreducible map. (d) The excised maps have outer degree $2b$ and girth $2b$. \label{fig:substitution}}
\end{figure}

In a rooted map the face to the right of the root edge is called the \emph{outer face} and all other faces are \emph{inner faces}. 
The \emph{outer degree} is the degree of its outer face.
Let $\vec{\mathcal{M}}_{0,\ell}^{(b)}$ be the set of $2b$-irreducible rooted planar maps with outer degree $2\ell$ with the additional requirement in the case $\ell=b$ that any cycle of length $2b$ bounds an inner face of degree $2b$.
Let also $\vec{\mathcal{G}}_{0,\ell}^{(\geq b)}$ be the set of rooted planar maps with outer degree $2\ell$ and girth at least $2b$.
Contrary to $\mathcal{M}_{0}^{(b)}$, we do include maps with one or two faces in $\vec{\mathcal{M}}_{0,\ell}^{(b)}$ and $\vec{\mathcal{G}}_{0,\ell}^{(\geq b)}$.
The corresponding generating functions are denoted
\begin{align*}
F^{(b)}_{0,\ell}(x_b,x_{b+1},\ldots) &= \sum_{\map\in\vec{\mathcal{M}}^{(b)}_{0,\ell}}\,\,\prod_{f\in \mathcal{F}'(\map)} x_{\deg(f)/2},\\
G^{(\geq b)}_{0,\ell}(x_b,x_{b+1},\ldots) &= \sum_{\map\in\vec{\mathcal{G}}^{(\geq b)}_{0,\ell}} \prod_{f\in \mathcal{F}'(\map)} x_{\deg(f)/2},
\end{align*}
where $\mathcal{F}'(\map)$ denotes the faces of $\map$ excluding the degree-$2\ell$ outer face.
The only $2b$-irreducible maps of outer degree $\ell \leq b$ are rooted plane trees and, in the case $\ell=b$, the rooted map consisting of a single cycle of length $2b$.
Hence
\begin{equation}\label{eq:Fconditions}
F_{0,\ell}^{(b)}(x_{b},x_{b+1}, \ldots) = \operatorname{Cat}(\ell) + x_b\, \ind_{\ell=b} \qquad \text{for }1\leq\ell\leq b. 
\end{equation}
It is shown in \cite[Section 3]{Bouttier2014} that these generating functions satisfy the relations
\begin{align}
G_{0,\ell}^{(\geq b)}(x_b,x_{b+1},\ldots) &= F_{0,\ell}^{(b-1)}(0,x_b,x_{b+1},\ldots), \label{eq:substitutionrel1}\\
G_{0,\ell}^{(\geq b)}(x_b,x_{b+1},\ldots) &= F_{0,\ell}^{(b)}(G_{0,b}^{(b)}(x_b,x_{b+1},\ldots),x_{b+1},x_{b+2},\ldots),\label{eq:substitutionrel2}\\
G_{0,b}^{(b)}(x_b,x_{b+1},\ldots) &= G_{0,b}^{(\geq b)}(x_b,x_{b+1},\ldots) -  \operatorname{Cat}(b).\label{eq:substitutionrel3}
\end{align}
The first identity expresses that the planar maps of girth at least $2b$ are precisely those maps that are $(2b-2)$-irreducible and have no inner face of degree $2b-2$.
The second and most non-trivial identity holds because any planar map of girth at least $2b$ can be obtained from a $2b$-irreducible planar map by gluing maps of girth $2b$ inside its faces of degree $2b$.
The last identity follows from the fact that a planar map of girth at least $2b$ and outer degree $2b$ has girth exactly $2b$ unless it is a tree.
These relations are illustrated in Figure~\ref{fig:substitution} in the case $b=2$.

According to \cite[Section 3.2]{Bouttier2014} the power series $x_b \mapsto G^{(b)}_{0,b}(x_b,x_{b+1},\ldots)$ has a compositional inverse, allowing one to determine $F_{0,\ell}^{(b)}$ in terms of $F_{0,\ell}^{(b-1)}$ for any $b\geq 1$.
By iteration one can find formal power series $X_j^{(b)}(x_b,x_{b+1}, \ldots)$, $j=1,\ldots,b$, such that
\begin{equation}\label{eq:planarsubstitution}
F_{0,\ell}^{(b)}(x_b,x_{b+1},\ldots) = F_{0,\ell}(X_1^{(b)},X_2^{(b)},\ldots,X_b^{(b)},x_{b+1}, \ldots),
\end{equation}  
where $F_{0,\ell}(x_1,x_2,\ldots) = F_{0,\ell}^{(0)}(0,x_1,x_2,\ldots)$ is the generating function of all (bipartite) planar maps with outer degree $2\ell$.

\subsection{Partition function of $2b$-irreducible planar maps}
Using universal properties of the generating function of planar maps, it possible to determine the left-hand side of \eqref{eq:planarsubstitution} without explicit constructing the power series $X_1^{(b)},\ldots,X_b^{(b)}$.
In fact, it is sufficient to know that such power series exist and that $F_{0,\ell}^{(b)}$ satisfies the conditions \eqref{eq:Fconditions}.
The results of \cite[Section 3.3]{Bouttier2014} may be summarized as follows.
Let 
\begin{equation}\label{eq:Uk}
U_k(x_1,x_2,\ldots) = \sum_{j=k+1}^\infty \binom{2j-1}{j+k}x_{j} R^{j+k},
\end{equation}
where $R(x_1,x_2,\ldots) = 1 + \cdots$ is the formal power series determined by
\begin{equation}
R = 1 + \sum_{j=1}^\infty \binom{2j-1}{j} x_{j}R^j.
\end{equation}
Then $F_{0,\ell}$ is given by
\begin{equation}
F_{0,\ell}(x_1,x_2,\ldots) = R^\ell \sum_{k=0}^\ell \frac{2k+1}{2\ell+1}\binom{2\ell+1}{\ell-k} (\ind_{\{k=0\}} - U_k R^{-k}\ind_{\{k>0\}}).
\end{equation}
When $b \geq 1$, one may perform the substitutions \eqref{eq:planarsubstitution} and introduce the power series
\begin{equation*}
R^{(b)}(x_b,x_{b+1},\ldots) \equiv R(X_1^{(b)},\ldots,X_{b}^{(b)},x_{b+1},\ldots), \quad U_k^{(b)}(x_b,x_{b+1},\ldots) \equiv U_k(X_1^{(b)},\ldots,X_{b}^{(b)},x_{b+1},\ldots).
\end{equation*}
The conditions \eqref{eq:Fconditions} for $0\leq \ell < b$ can be solved for $U_k^{(b)}$, leading to the expressions \cite[Equations (3.15) \& (3.17)]{Bouttier2014}
\begin{align}
U_k^{(b)} &=\begin{cases}
-\sum_{m=0}^k\, (-1)^{k+m}\binom{k+m}{2m}\operatorname{Cat}(m) (R^{(b)})^{k-m}&\text{for }1\leq k < b,\\
\sum_{j=k+1}^\infty\binom{2j-1}{j+k}x_{j}(R^{(b)})^{j+k}&\text{for }k\geq b.
\end{cases}\label{eq:Usubstitution}
\end{align}
Finally the condition \eqref{eq:Fconditions} for $\ell=b$ uniquely determines $R^{(b)}$ by
\begin{align}
Z^{(b)}(R^{(b)};x_{b+1},\cdots) &= x_b, \nonumber\\
Z^{(b)}(r;x_{b+1},\cdots) &\coloneqq r\,\ind_{\{b=0\}}-\sum_{\ell=0}^b (-1)^{b-\ell}\binom{b+\ell}{2\ell}\operatorname{Cat}(\ell)r^{b-\ell} - \sum_{j=b+1}^\infty \binom{2j-1}{j+b}x_{j}r^{b+j}.\label{eq:Zfun}
\end{align}

The first term of $Z^{(b)}$ is included to obtain a correct criterion in the case $b=0$ as well.
To see this, note that 
\begin{equation}\label{eq:Z0equation}
x_0 = Z^{(0)}(R^{(0)};x_1,\ldots) = R^{(0)} - 1 - \sum_{j=1}^\infty \binom{2j-1}{j}x_{j}(R^{(0)})^{j}
\end{equation} 
is solved by
\begin{equation}\label{eq:R0}
R^{(0)}(x_0,x_1,\ldots) = (1+x_0)\, R\left( x_1,\,(1+x_0)\,x_2,\,(1+x_0)^2 \,x_3,\ldots\right).
\end{equation}
By Euler's formula 
\begin{equation*}
|\mathcal{V}(\map)| = 2 + |\mathcal{E}(\map)| - |\mathcal{F}(\map)| = 2 + \sum_{f\in\mathcal{F}(\map)} (\tfrac{1}{2}\deg(f)-1),
\end{equation*}
the substitution $x_k \mapsto (1+x_0)^{k-1} x_k$ for $k\geq 1$ is the correct one to insert a weight $(1+x_0)$ per vertex.
We will see in the Proposition below that the overall normalization of the partition function also comes out right in this case.

Combining various results of \cite{Bouttier2014} leads to an explicit expression for the second derivative of $F_0^{(b)}$, i.e.\ the partition function with two distinguished faces.

\begin{proposition}\label{thm:genus0partitionfunction}
For $b\geq 0$ the partition function of $2b$-irreducible planar maps with at least three faces satisfies for $\ell,\ell'\geq \max(b,1)$,
\begin{equation}\label{eq:genus0partitionfunction}
\frac{\partial^2 F_0^{(b)}}{\partial x_{\ell_1}\partial x_{\ell_2}}(x_b,x_{b+1},\ldots) = \frac{3x_b^2+2x_b^3}{2(1+x_b)^2}\,\ind_{\ell_1=\ell_2=b}  +  \int_{1+x_0\ind_{\{b=0\}}}^{R^{(b)}} \rmd r \frac{\prod_{i=1}^2 \left(\ind_{\ell_i=b} + \binom{2\ell_i-1}{\ell_i+b}r^{\ell_i+b}\right)}{r^{2b+1}}.
\end{equation}
\end{proposition}
\begin{proof}
We treat four cases separately: $\ell_1,\ell_2 > b>0$; $\ell_1 > \ell_2= b>0$; $\ell_1=\ell_2=b>0$; $\ell_1, \ell_2 > b=0$.

\noindent\textbf{Case $\ell_1,\ell_2>b$:} According to \cite[Equation (9.17)]{Bouttier2014} the generating function of $2b$-irreducible maps with two distinguished faces of degrees $2\ell_1,2\ell_2 > 2b$ with a marked edge on both faces reads
\begin{equation*}
(\ell_1-b)(\ell_2-b) \binom{2\ell_1}{\ell_1-b}\binom{2\ell_2}{\ell_2-b}\frac{(R^{(b)})^{\ell_1+\ell_2}}{\ell_1+\ell_2}.
\end{equation*}
Subtracting the contribution of the maps with only two faces ($R^{(b)}=1$) and compensating for the marked edges we thus have
\begin{align}
\frac{\partial^2 F_0^{(b)}}{\partial x_{\ell_1}\partial x_{\ell_2}} &= \frac{1}{2\ell_1}\frac{1}{2\ell_2} (\ell_1-b)(\ell_2-b) \binom{2\ell_1}{\ell_1-b}\binom{2\ell_2}{\ell_2-b}\frac{(R^{(b)})^{\ell_1+\ell_2}-1}{\ell_1+\ell_2}\nonumber\\
&=
\binom{2\ell_1-1}{\ell_1+b}\binom{2\ell_2-1}{\ell_2+b}\frac{(R^{(b)})^{\ell_1+\ell_2}-1}{\ell_1+\ell_2},\label{eq:cylindergenfun}
\end{align}
which agrees with the right-hand side of \eqref{eq:genus0partitionfunction}.

\noindent\textbf{Case $\ell_1>\ell_2=b$:} According to \cite[Equation (3.26)]{Bouttier2014},
\begin{equation*}
\frac{\partial F^{(b)}_{0,\ell_1}}{\partial x_b} = \binom{2\ell_1}{\ell_1-b} (R^{(b)})^{\ell_1-b}.
\end{equation*} 
Subtraction of the contribution with just two faces and compensation for the root on the outer face yields 
\begin{equation*}
\frac{\partial^2 F^{(b)}_0}{\partial x_{\ell_1}\partial x_{\ell_2}} = \frac{1}{2\ell_1}  \binom{2\ell_1}{\ell_1-b} \left((R^{(b)})^{\ell_1-b}-1\right) = \binom{2\ell_1-1}{\ell_1+b} \frac{(R^{(b)})^{\ell_1-b}-1}{\ell_1-b},
\end{equation*}
in agreement with \eqref{eq:genus0partitionfunction}.

\noindent\textbf{Case $\ell_1=\ell_2=b$:} According to \cite[Equation (9.1)]{Bouttier2014} the generating function of rooted $2b$-irreducible maps with outer face of degree $2b$ is given by
\begin{equation*}
H_b(x_b,x_{b+1},\ldots) = 2x_b + \frac{b x_b^3}{1+x_b} - X_b(x_b,x_b+1,\ldots),
\end{equation*}
where $x_b\mapsto X_b(x_b,x_{b+1},\ldots)$ is the functional inverse of $x_b\mapsto G_{0,b}^{(b)}(x_b,x_{b+1},\ldots)$.
With the help of \eqref{eq:substitutionrel1}, \eqref{eq:substitutionrel3} and \eqref{eq:cylindergenfun} we find
\begin{equation*}
\frac{\partial G_{0,b}^{(b)}}{\partial x_b} = \frac{\partial F_{0,b}^{(b-1)}}{\partial x_b} = (R^{(b-1)})^{2b}.
\end{equation*}
It follows that
\begin{equation*}
\frac{\partial X_b}{\partial x_b} (x_b,x_{b+1},\ldots) = \left(\frac{\partial G_{0,b}^{(b)}}{\partial x_b}(X_b,x_{b+1},\ldots)\right)^{-1} =  \left(R^{(b-1)}(X_b,x_{b+1},\ldots)\right)^{-2b} = \left(R^{(b)}(x_b,x_{b+1},\ldots)\right)^{-2b},
\end{equation*}
where the last equality is a direct consequence of the substitution approach.
Hence, 
\begin{align*}
\frac{\partial^2 F^{(b)}_0}{\partial x_b^2} = \frac{1}{2b} \frac{\partial H_b}{\partial x_b} - \frac{1}{2b} \frac{\partial H_b}{\partial x_b}(0,0,\ldots) = \frac{3x_b^2+2x_b^3}{2(1+x_b)^2} + \frac{1-(R^{(b)})^{-2b}}{2b}, 
\end{align*}
again agrees with \eqref{eq:genus0partitionfunction}.

\noindent\textbf{Case $\ell_1,\ell_2>b=0$:} It is well-known (see for example \cite[Theorem 1.1]{Collet2012}) that the generating function of arbitrary bipartite planar maps with even face degrees and two distinguished faces of degree $2\ell_1,2\ell_2 \geq 2$ and weight $1+x_0$ per vertex is given by
\begin{equation*}
\binom{2\ell_1-1}{\ell_1}\binom{2\ell_1-1}{\ell_1} \frac{(R^{(0)})^{\ell_1+\ell_2}}{\ell_1+\ell_2},
\end{equation*}
where $R^{(0)}$ satisfies \eqref{eq:Z0equation}.
Subtracting the contribution with only two faces ($R^{(0)}=1+x_0$), this agrees with \eqref{eq:genus0partitionfunction}.
\end{proof}

Apart from a small correction when all faces are of degree $2b$ the partition function $\frac{\partial^2 F_0}{\partial x_{\ell_1}\partial x_{\ell_2}}$ with two distinguished faces of degree $2\ell_1$ and $2\ell_2$ is completely expressed in terms of 
\begin{equation}\label{eq:Rinterpretation}
R^{(b)}(x_b,x_{b+1},\ldots) = 1+\frac{\partial^2 F_0^{(b)}}{\partial x_{b+1}\partial x_{b}}(x_b,x_{b+1},\ldots).
\end{equation}
This enumerates all $2b$-irreducible maps with two marked faces of degrees $2b+2$ and $2b$, where the $1$ takes into account the unique such map with only two faces.
In \cite[Section 4]{Bouttier2014} the quantity $R^{(b)}- 1$ is interpreted as the generating function of certain $2b$-irreducible \emph{slices} and a combinatorial interpretation of its equation \eqref{eq:Zfun} can be understood via a \emph{slice decomposition}.
This interpretation will not play a role in this work, but it is certainly something one would like to understand better for other topologies and when vertices of degree one are disallowed.
	
\subsection{Substitution approach for higher genus}

Recall that we call $\map$ essentially $2b$-irreducible if $\map^\infty$ is $2b$-irreducible.
The generating function of such maps is given by $F_g^{(b)}$ in \eqref{eq:partitionfunction}.
Similarly we say $\map$ has \emph{essential girth} $2b$ if $\map^{\infty}$ has girth $2b$.
Denoting the set of rooted, respectively unrooted, maps of essential girth at least $2b$ by $\vec{\mathcal{G}}^{(\geq b)}_{g}$ respectively $\mathcal{G}^{(\geq b)}_{g}$, we introduce the generating function
\begin{align*}
G_g^{(\geq b)}(x_b,x_{b+1}, \ldots) &= \sum_{\map\in\mathcal{G}^{(b)}_g} \frac{1}{|\operatorname{Aut}(\map)|} \prod_{f\in \mathcal{F}(\map)} x_{\deg(f)/2}  = \sum_{\map\in\vec{\mathcal{G}}^{(b)}_g} \frac{1}{2|\mathcal{E}(\map)|} \prod_{f\in \mathcal{F}(\map)} x_{\deg(f)/2}.
\end{align*}
The substitution approach of Bouttier \& Guitter extends to the following relation between these generating functions.
\begin{proposition}
	For $g\geq 1$ and $b\geq 1$ we have the identities of formal power series
	\begin{equation}\label{eq:girthgenfunidentities}
	F_{g}^{(b-1)}(0,x_b,x_{b+1},\ldots) = G_{g}^{(\geq b)}(x_b,x_{b+1}, \ldots) = F_{g}^{(b)}(G^{(b)}_{0,b}(x_b,x_{b+1},\ldots),x_{b+1},x_{b+2},\ldots). 
	\end{equation}
\end{proposition}

\begin{proof}
The first identity just expresses the fact that the maps of essential girth at least $2b$ are precisely the essentially $(2b-2)$-irreducible maps that have no faces of degree $2b-2$.
To deduce the second identity we will show that these maps can also be obtained by gluing maps inside the faces of degree $2b$ of essentially $2b$-irreducible maps.
Denoting by $\mathcal{F}_b(\map)$ the set of faces of $\map$ of degree $2b$, let us consider the mapping 
\begin{equation*}
\mathsf{Glue} : \Big\{ (\map, (\map_f)_{f \in \mathcal{F}_b(\map)}) : \map\in \vec{\mathcal{M}}^{(b)}_{g},\, \map_f \in \vec{\mathcal{G}}^{(b)}_{0,b}\Big\} \to  \vec{\mathcal{G}}^{(\geq b)}_{g}
\end{equation*}
defined by taking $\map' = \mathsf{Glue}( \map, (\map_f)_{f \in \mathcal{F}_b(\map)})$ to be the rooted map obtained from $\map$ by gluing $\map_f$ inside the face $f$ of $\map$ for each face $f$ of degree $2b$.
To perform the gluing we select an arbitrary but deterministic algorithm to distinguish an edge on each face $f$ of degree $2b$ to which the root edge of $\map_f$ is glued. 
Let ${\map'}^\infty$ be the universal cover of $\map'$.
Applying the encircling lemma of \cite[Section 3.1]{Bouttier2014} to ${\map'}^\infty$ shows that ${\map'}^\infty$ has girth at least $2b$, so indeed $\map' \in  \vec{\mathcal{G}}^{(\geq b)}_{g}$.

Now let $\map'\in \vec{\mathcal{G}}^{(\geq b)}_{g}$ be an arbitrary rooted map of essential girth at least $2b$ and ${\map'}^\infty$ its universal cover.
As in Lemma \ref{lem:irrcriterion}, let $\mathcal{C}_d({\map'}^\infty)$ be the set of \emph{outermost} cycles of length $2b$ in ${\map'}^\infty$.
If the essential girth of $\map'$ is larger than $2b$ then $\mathcal{C}_d({\map'}^\infty)$ is empty, while if the essential girth is $2b$ the proof of Lemma~\ref{lem:irrcriterion} shows that $\mathcal{C}_d({\map'}^\infty)$ comprises of the lifts of a finite collection of essentially simple cycles $c_1,\ldots,c_p$ in $\map'$ that have non-overlapping interiors.
If the root edge of $\map'$ is not contained in one of these interiors, then we may unambiguously construct a rooted map $\map''$ by removing the interiors of these cycles, leaving faces of degree $2b$ in their place.
By construction, its universal cover ${\map''}^\infty$ has girth $2b$ and the only cycles of length $2b$ are the boundaries of a face of degree $2b$, so $\map''$ is essentially $2b$-irreducible.
The regions that have been excised from $\map'$ are by construction planar, of girth $2b$, and have outer degree $2b$.

If $\map' = \mathsf{Glue}( \map, (\map_f)_{f \in \mathcal{F}_d(\map)})$ then the cycles $c_1,\ldots,c_p$ are easily seen to correspond to the boundaries of the maps $\map_f$ inside $\map'$.
Hence the construction above reproduces the original map, i.e. $\map''=\map$, while the excised maps are precisely $(\map_f)_{f \in \mathcal{F}_d(\map)}$.
Therefore $\mathsf{Glue}$ is injective and the image of $\mathsf{Glue}$ is given by the subset $\vec{\mathcal{G}}^{(\geq b)}_{g,\mathrm{out}}\subset\vec{\mathcal{G}}^{(\geq b)}_{g}$ of maps whose root edge is not contained in the interior of a cycle of length $2b$.

As explained in \cite[Section 3.1]{Bouttier2014} the mapping $\mathsf{Glue}$ preserves the collection of face degrees via the identities
\begin{equation*}
|\mathcal{F}_b(\map')| = \sum_{f\in\mathcal{F}_b(\map)} |\mathcal{F}'_b(\map_f)|, \qquad |\mathcal{F}_k(\map')| = |\mathcal{F}_k(\map)| + \sum_{f\in\mathcal{F}_b(\map)} |\mathcal{F}'_k(\map_f)|\quad \text{for }k>b,
\end{equation*}
where $\mathcal{F}_k(\map)$ denotes the set of faces of $\map$ of degree $2k$ and similarly for $\mathcal{F}'_k(\map)$ but excluding the outer face.
It follows that we have
\begin{align*}
F_{g}^{(b)}(G^{(b)}_{0,b}(x_d,x_{d+1},\ldots),x_{d+1},x_{d+2},\ldots) & = \!\!\!\!
\sum_{\substack{\map\in\vec{\mathcal{M}}^{(b)}_g\\(m_f)_{f\in \mathcal{F}_b(m)}}}\!\! \frac{1}{2|\mathcal{E}(\map)|} x_b^{\sum_{f\in\mathcal{F}_b(\map)} |\mathcal{F}'_b(\map_f)|}  \prod_{j>b}x_j^{|\mathcal{F}_j(\map)| + \sum_{f\in\mathcal{F}_b(\map)} |\mathcal{F}'_j(\map_f)|}\\
&= \sum_{\map'\in \vec{\mathcal{G}}^{(\geq b)}_{g,\mathrm{out}}}\frac{1}{2|\mathcal{E}_{\mathrm{out}}(\map')|}\prod_{j\geq b}x_j^{|\mathcal{F}_j(\map')|}\\
&= \sum_{\map'\in \vec{\mathcal{G}}^{(\geq b)}_{g}}\frac{1}{2|\mathcal{E}(\map')|}\prod_{j\geq b}x_j^{|\mathcal{F}_j(\map')|}\\
&= G_{g}^{(\geq b)}(x_b,x_{b+1}, \ldots),
\end{align*}
where $\mathcal{E}_{\mathrm{out}}(\map')$ denotes the set of edges of $\map'$ that are not contained in the interior of a cycle of length $2b$.
This proves the second identity in \eqref{eq:girthgenfunidentities}.
\end{proof}

As in the planar case one may use the compositional inverse of the power series $x_b \mapsto G^{(b)}_{0,d}(x_b,x_{b+1},\ldots)$ to determine $F_{g}^{(b)}$ in terms of $F_{g}^{(b-1)}$ for any $b\geq 1$.
By iteration one finds that
\begin{equation}\label{eq:highergenussubtitution}
F_{g}^{(b)}(x_b,x_{b+1}, \ldots) = F^{(0)}_{g}(0,X_1^{(b)},X_2^{(b)},\ldots,X_b^{(b)},x_{b+1}, \ldots) =  F_{g}(X_1^{(b)},X_2^{(b)},\ldots,X_b^{(b)},x_{b+1}, \ldots),
\end{equation}  
where $X_j^{(b)}(x_b,x_{b+1},\ldots)$, $1\leq j \leq b$, are the same formal power series as appeared in \eqref{eq:planarsubstitution} in the planar case.

\subsection{Partition function for $2b$-irreducible genus-$g$ maps}\label{sec:highergenuspartition}
The partition function $F_g(x_1,x_2,\ldots)$ of genus-$g$ maps with even faces can be extracted from topological recursion, see \cite[Chapter 3]{Eynard2016} for an overview.
The relevant spectral curve is given in terms of $U_k(x_1,x_2,\ldots)$ from \eqref{eq:Uk} and $R(x_1,x_2,\ldots)$ by\footnote{Note that \cite[Chapter 3]{Eynard2016} deals more generally with maps that are not necessarily bipartite and that the notation differs from the one used here, although the curves agree. For convenience we include a dictionary (where the left-hand sides are Eynard's notation): 
\begin{align*}
	&t=1, \qquad t_{2j} = x_{j}, \qquad \gamma = \sqrt{R},\qquad  u_1 = R^{-1/2}, \qquad M_{+,0} = \bar{M}_0/R, \\ 
	&u_{2k+1} = - U_k R^{-k-1/2}, \qquad u_{2k}=0, \qquad M_{+,k} = (-1)^k M_{-,k} = - R^{-k/2} \bar{M}_k / \bar{M}_0 \qquad \text{ for }k\geq 1.
\end{align*}
}
\begin{equation}
\begin{cases}
x(z) = \sqrt{R}(z+\tfrac{1}{z})\\
y(z)= -\frac{1}{2\sqrt{R}}\left(z-\tfrac{1}{z}-\sum_{k=1}^\infty U_k R^{-k} (z^{2k+1}-z^{-2k-1})\right).
\end{cases}
\end{equation}
Alternatively one may express $y(z)$ as
\begin{equation}
y(z) = -\frac{1}{2\sqrt{R}}\,(z-\tfrac{1}{z})\sum_{p=0}^\infty \bar{M}_p\,(z+\tfrac{1}{z}-2)^p
\end{equation}
in terms of the \emph{moments}
\begin{equation}\label{eq:moments}
\bar{M}_p(x_1,x_2,\ldots) = \ind_{p=0} - \sum_{k=0}^\infty \binom{2k+p+1}{2p+1} U_k\, R^{-k}.
\end{equation}
Note that for $p\geq 0$ fixed, $\binom{2k+p+1}{2p+1}$ may be interpreted as a polynomial in $k$ of degree $2p+1$ that vanishes at all integer values with absolute value smaller than $p/2$.
According to \cite[Theorem 3.4.6 \& Corollary 3.5.1]{Eynard2016} the partition functions $F_g$ for $g\geq 1$ are expressible in terms of the moments via
\begin{equation}\label{eq:Fgmaps}
F_1 = - \frac{1}{12} \log \bar{M}_0, \qquad F_g = P_g\left( \frac{1}{\bar{M}_0}, \frac{\bar{M}_1}{\bar{M}_0},\ldots, \frac{\bar{M}_{3g-3}}{\bar{M}_0}\right) \quad \text{for }g\geq 2,
\end{equation}
where $P_g(m_0,\ldots,m_{3g-3})$ is a universal polynomial with rational coefficients for each $g\geq 2$.
More precisely $P_g$ is of the form 
\begin{equation}\label{eq:Ptilde}
P_g(m_0,\ldots,m_{3g-3}) = \tilde{P}_g(0,\ldots,0) - m_0^{2g-2} \tilde{P}_g(m_1,\ldots,m_{3g-3})
\end{equation}
for some polynomial $\tilde{P}_g(m_1,\ldots,m_{3g-3})$ such that $\tilde{P}_g(\mu,\mu^2,\ldots,\mu^{3g-3})$ is of degree $3g-3$ in $\mu$.
For example, in the genus-$2$ case we find with the help of \cite[Section 3.5.2]{Eynard2016} the polynomial
\begin{equation*}
P_2(m_0,m_1,m_2,m_3) = \frac{1}{240}-\frac{m_0^2}{30720} \left(2016 m_1^3+1086 m_1^2-3480 m_2 m_1+407 m_1-860 m_2+1400 m_3+128\right).
\end{equation*}
Similar results have been obtained from Hermitian matrix models, see  \cite{Ambjorn1993,Akemann1996}.

The expression \eqref{eq:Fgmaps} in terms of moments is convenient because we know the result of substituting $x_{j} \to X^{(b)}_{j}(x_b,x_{b+1},\ldots)$ for $1<j \leq b$ in $R$ and $U_k$.
With the help of \eqref{eq:Usubstitution} and \eqref{eq:moments} we find for $b\geq 1$ that
\begin{align}
\bar{M}_p^{(b)}(x_b,x_{b+1},\ldots) &\equiv \bar{M}_p(X_1^{(b)},\ldots,X_b^{(b)},x_{b+1},\ldots)\nonumber\\
&=\ind_{p=0} - \sum_{k=1}^\infty \binom{2k+p+1}{2p+1} U_k^{(b)}\, (R^{(b)})^{-k} \label{eq:Mbarbdef}\\
&= \sum_{k=0}^{b-1} \binom{2k+p+1}{2p+1} \sum_{m=0}^k (-1)^{k+m}\binom{k+m}{2m}\operatorname{Cat}(m) (R^{(b)})^{-m} \nonumber\\
&\quad - \sum_{k=b}^\infty \binom{2k+p+1}{2p+1}\sum_{j=k+1}^\infty \binom{2j-1}{j+k}x_{j}(R^{(b)})^{j}.\label{eq:mbarexpr}
\end{align}
To rewrite these sums we make use of the following lemma.

\begin{lemma}\label{lem:Qsums}
For each $p\geq 0$, there exists a unique polynomial $Q_p(b,j)$ such that
\begin{align}\label{eq:sumpolynomial}
\sum_{k=b}^{j-1} \binom{2k+1+p}{2p+1}\binom{2j-1}{j+k} &= \binom{2j-1}{j+b} Q_p(b,j) \qquad \text{for }j > b \geq 0, \\
\sum_{k=m}^{b-1} \binom{2k+1+p}{2p+1}(-1)^{k+m}\binom{k+m}{2m} &= -(-1)^{b+m}\binom{b+m}{2m} Q_p(b,-m)\qquad \text{for }b > m \geq 0.\label{eq:sumpolynomial2}
\end{align}
The polynomial $Q_p(b,j)$ is of degree $2p+1$ in $b$ and degree $p+1$ in $j$ and satisfies $Q_p(b,-b)=0$.
\end{lemma}
\begin{proof}
	We apply Gosper's algorithm \cite{Gosper1978} to the sum in \eqref{eq:sumpolynomial}. 
	Let us denote the summand by
	\begin{equation*}
	t(k) = \binom{2k+1+p}{2p+1}\binom{2j-1}{j+k}
	\end{equation*}
	and introduce the polynomials
	\begin{equation*}
	\quad p(k) = \binom{2k+1+p}{2p+1},\quad q(k) = j-k-1, \quad r(k) = j+k.
	\end{equation*}
	Then the ratio of consecutive summands satisfies
	\begin{equation}
	\frac{t(k+1)}{t(k)} = \frac{q(k)}{r(k+1)} \frac{p(k+1)}{p(k)}.
	\end{equation}
	According to Gosper, there exists a solution to $T(k+1)-T(k)=t(k)$ of the form of a hypergeometric term if and only if there exists a polynomial $s(k)$ of degree $2p$ that solves the recurrence equation
	\begin{equation}\label{eq:gosperrecur}
	p(k) = q(k)s(k+1)-r(k)s(k).
	\end{equation}

	First we observe that $s(k) \mapsto q(k) s(k+1) - r(k) s(k)$ determines an injective linear mapping (working in the ring of polynomials in $j$) from polynomials of degree $2p$ to those of degree $2p+1$. 
	We thus need to determine the one-dimensional cokernel of the mapping, in other words we need a single linear condition satisfied by all polynomials of the form $q(k) s(k+1) - r(k) s(k)$.
	A convenient way to achieve this is to change the polynomial to a differential operator and to seek a formal power series $V(x)$ satisfying
	\begin{align*}
	0 &= \left[s(1+\partial_x)q(\partial_x) - s(\partial_x)r(\partial_x) \right] V(x)\big|_{x=0} \\
	&= \left[ e^{-x} s(\partial_x) e^{x} ((j-1)-\partial_x)- s(\partial_x)(j+\partial_x)\right]  V(x)\big|_{x=0} \\
	&= s(\partial_x)\left[ e^{x} ((j-1)-\partial_x)-(j+\partial_x)\right] V(x)\big|_{x=0}.
	\end{align*}
	It is now easy to find a power series solution independent of $s(k)$,
	\begin{equation*}
	V(x) = e^{-x/2} (\cosh(x/2))^{2j-1},
	\end{equation*}
	which has coefficients that are polynomials in $j$.
	To determine whether $p(k)$ is in the image, we need to check that $p(\partial_x) V(x)|_{x=0} = 0$.
	We may calculate
	\begin{align*}
	p(\partial_x) V(x)|_{x=0} &= p(\partial_x-\tfrac{1}{2}) e^{x/2} V(x)|_{x=0}\\
	&= p(\partial_x-\tfrac{1}{2})(\cosh(x/2))^{2j-1}|_{x=0},
	\end{align*}
	but this vanishes because $p(\ell - \tfrac12)$ is an odd polynomial in $\ell$ while the formal power series $(\cosh(x/2))^{2j-1}$ is even in $x$.
	We conclude that there exists a polynomial $s(k)$ of degree $2p$ in $k$ whose coefficients are polynomials in $j$ solving \eqref{eq:gosperrecur}.
	
	In this case $T(k) = r(k) t(k) s(k) / p(k)$ and therefore \eqref{eq:sumpolynomial} is satisfied with
	\begin{equation}\label{eq:Qgosper}
	Q_p(b,j) = \frac{T(j)-T(b)}{\binom{2j-1}{j+b}} = -r(b) s(b) = -(j+b) s(b),
	\end{equation}
	which is indeed polynomial in $b$ and $j$. 	
	This finishes the proof of the first identity.
	It also follows directly from \eqref{eq:Qgosper} that $Q_p(b,j)$ is of degree $2p+1$ in $b$ and satisfies $Q_p(b,-b)=0$.
	
	In order to determine the degree of $Q_p(b,j)$ in $j$ we perform an asymptotic analysis of the summand in \eqref{eq:sumpolynomial} as $k,j\to\infty$ for $b,p$ fixed.
	Using Stirling's approximation we easily find that there exists a smooth convex function $f: (0,\infty) \to \R$ with global maximum at $\sqrt{p+1/2}$ such that
	\begin{align*}
	\log\left(\frac{\binom{2k+1+p}{2p+1}\binom{2j-1}{j+k}}{\binom{2j-1}{j+b} \,j^{p+\tfrac{1}{2}}}\right) \xrightarrow[k\sim y\sqrt{j}]{j,k\to\infty} f(y).
	\end{align*}
	Approximating the sum in \eqref{eq:sumpolynomial} by an integral, this implies that $\log Q_p(b,j) \sim (p+1) \log j$ as $j\to\infty$ for any $p,b$ fixed and therefore $Q_p(b,j)$ is of degree $p+1$ in $j$.
	
	For the second sum, note that 
	\begin{align*}
	\binom{2k+1+p}{2p+1}(j+k+1) &= p(k) r(k+1) \stackrel{\eqref{eq:gosperrecur}}{=}q(k)s(k+1)r(k+1) - r(k+1)r(k)s(k) \\
	&\stackrel{\mathclap{\eqref{eq:Qgosper}}}{=} r(k+1) Q_p(k,j) - q(k) Q_p(k+1,j)\\
	&= (j+k+1) Q_p(k,j) - (j-k-1) Q_p(k+1,j)
	\end{align*}
	is an equality between polynomials in $j$ and $k$.
	Setting $j=-m$, we easily find that 
	\begin{align*}
	\binom{2k+1+p}{2p+1} (-1)^{k+m}\binom{k+m}{2m} = (-1)^{k+m}\binom{k+m}{2m} Q_p(k,-m) - (-1)^{k+m+1}\binom{k+1+m}{2m} Q_p(k+1,-m)
	\end{align*}
	for $k\geq m\geq 0$.
	Summing over $k$ from $m$ to $b-1$ and using that $Q_p(m,-m)=0$ yields precisely the identity \eqref{eq:sumpolynomial2}.
\end{proof}

\noindent
The first few polynomials are given by
\begin{align*}
Q_0(b,j)&=b+j,\\
Q_1(b,j)&=\frac{2}{3} (b+j) \left(b^2+j-1\right),\\
Q_2(b,j)&=\frac{1}{30} (b+j) \left(4 b^4+b^2 (8 j-15)+(j-1) (8 j-11)\right),\\
Q_3(b,j)&=\frac{1}{315} (b+j) \left(4 b^6+b^4 (12 j-35)+b^2 (24 j^2-90 j+91)+6
(j-1)(j-2)  (4 j-5)\right).
\end{align*}
With these polynomials in hand the moments can be expressed rather concisely.

\begin{proposition}\label{thm:highergenuspartitionfunctions}
The partition functions $F^{(b)}_g(x_b,x_{b+1},\ldots)$ for $g\geq 1$ and $b\geq 0$ are given by
\begin{align}
F^{(b)}_1 &= - \frac{1}{12} \log \bar{M}^{(b)}_0, \qquad F^{(b)}_g = P_g\left( \frac{1}{\bar{M}^{(b)}_0}, \frac{\bar{M}^{(b)}_1}{\bar{M}^{(b)}_0},\ldots, \frac{\bar{M}^{(b)}_{3g-3}}{\bar{M}^{(b)}_0}\right) \quad \text{for }g\geq 2,\,b\geq 1,\\
F^{(0)}_1 &= - \frac{1}{12} \log \frac{\bar{M}^{(0)}_0}{1+x_0}, \qquad F^{(0)}_g = (1+x_0)^{2-2g}P_g\left( \frac{1}{\bar{M}^{(0)}_0}, \frac{\bar{M}^{(0)}_1}{\bar{M}^{(0)}_0},\ldots, \frac{\bar{M}^{(0)}_{3g-3}}{\bar{M}^{(0)}_0}\right) \quad \text{for }g\geq 2,\label{eq:Fgbzero}
\end{align}
where
\begin{equation}\label{eq:mbarb}
\bar{M}_p^{(b)}(x_b,x_{b+1},\ldots) = Q_p(b,r\, \partial_r)\,r^{-b} \, Z^{(b)}(r; x_{b+1}, x_{b+2}, \ldots) \Big|_{r=R^{(b)}(x_b,x_{b+1},\ldots)}
\end{equation}
with $Z^{(b)}$ as given in \eqref{eq:Zfun}. 
\end{proposition}
\begin{proof}
With the help of Lemma \ref{lem:Qsums} the expression \eqref{eq:mbarexpr} for $b\geq 1$ evaluates to
\begin{align*}
\bar{M}_p^{(b)}(x_b,x_{b+1},\ldots) &= \sum_{m=0}^{b-1} \sum_{k=m}^{b-1}\binom{2k+p+1}{2p+1}  (-1)^{k+m}\binom{k+m}{2m}\operatorname{Cat}(m) (R^{(b)})^{-m} \nonumber\\
&\quad - \sum_{j=b+1}^\infty \sum_{k=b}^{j-1} \binom{2k+p+1}{2p+1} \binom{2j-1}{j+k}x_{2j}(R^{(b)})^{j}.\\
&= - \sum_{m=0}^{b-1} Q_p(b,-m) (-1)^{b+m}\binom{b+m}{2m}\operatorname{Cat}(m) (R^{(b)})^{-m} \\
&\quad - \sum_{j=b+1}^\infty Q_p(b,j) \binom{2j-1}{j+b}x_{2j}(R^{(b)})^{j}.\\
&= Q_p(b,r\, \partial_r) \,r^{-b}\, Z^{(b)}(r; x_{b+1}, x_{b+2}, \ldots) \Big|_{r=R^{(b)}(x_b,x_{b+1},\ldots)}.
\end{align*}
Performing the substitution \eqref{eq:highergenussubtitution} in \eqref{eq:Fgmaps} gives the desired result for $b \geq 1$.

It remains to check the case $b=0$. 
Euler's formula
\begin{equation*}
|\mathcal{V}(\map)| = 2 - 2g + |\mathcal{E}(\map)| - |\mathcal{F}(\map)| = 2-2g + \sum_{f\in\mathcal{F}(\map)} (\tfrac{1}{2}\deg(f)-1)
\end{equation*}
implies that
\begin{equation*}
F_g^{(0)}(x_0,x_1,\ldots) = (1+x_0)^{2-2g}F_g(x_1,\,(1+x_0)x_2,\,(1+x_0)^2x_3,\ldots).
\end{equation*}
Defining $\bar{M}_p^{(0)}$ in terms of $\bar{M}_p$ via 
\begin{equation}
\bar{M}_p^{(0)}(x_0,x_1,\ldots) = (1+x_0) \,\bar{M}_p(x_1,\,(1+x_0)x_2,\,(1+x_0)^2x_3,\ldots),
\end{equation}
we see that $F_g^{(0)}$ is indeed given by \eqref{eq:Fgbzero}. 
To check that $\bar{M}_p^{(0)}$ is given by \eqref{eq:mbarb}, we note that \eqref{eq:Usubstitution} and \eqref{eq:R0} together imply that
\begin{equation*}
U_k^{(0)}(x_0,x_1,\ldots) = (1+x_0)^{k+1} U_k(x_1,\,(1+x_0)x_2,\,(1+x_0)^2x_3,\ldots).
\end{equation*}
Therefore 
\begin{align*}
\bar{M}_p^{(0)}(x_0,x_1,\ldots) &=(1+x_0)\left( \ind_{\{p=0\}} - \sum_{k=1}^\infty \binom{2k+1+p}{2p+1} U_k(x_1,\,(1+x_0)x_2,\ldots)\,R(x_1,\,(1+x_0)x_2,\ldots)^{-k}\right)\\
&=R^{(0)}\,\ind_{\{p=0\}} - \sum_{k=0}^\infty \binom{2k+1+p}{2p+1} U_k^{(0)}\,(R^{(0)})^{-k},
\end{align*}
where we used that $U_0^{(0)} = R^{(0)} - (1+x_0)$ by \eqref{eq:Z0equation}.
Now we may insert the expression \eqref{eq:Usubstitution} for $U^{(0)}_k$ and apply Lemma \ref{lem:Qsums}. 
Since $Q_p(0,1) = \ind_{\{p=0\}}$, in this way we obtain
\begin{align*}
\bar{M}_p^{(0)}(x_0,x_1,\ldots)&=Q_p(0,r\,\partial_r)\left(r-1 - \sum_{j=1}^\infty \binom{2j-1}{j}x_j r^j\right) \Big|_{r=R^{(0)}(x_0,x_{1},\ldots)}\\
&\stackrel{\eqref{eq:Z0equation}}{=}Q_p(0,r\,\partial_r) Z^{(0)}(r;x_1,x_2,\ldots) \Big|_{r=R^{(0)}(x_0,x_{1},\ldots)}.
\end{align*}
This concludes the proof.
\end{proof}

\noindent
One may even express the moments entirely in terms of $R^{(b)}$ and its $x_b$-derivatives as is shown in the next lemma.

\begin{lemma}\label{lem:momentRder}
	For any $p\geq 0$, there exists a polynomial $T_p(b,r_0,\cdots,r_{p+1})$ homogeneous of degree $2p$ in $r_0,\ldots,r_{p+1}$ such that for any $b\geq 0$,
	\begin{equation}
	\bar{M}^{(b)}_p(x_{b},x_{b+1},\ldots) = \frac{(R^{(b)})^{1-b}}{(\partial_{x_{b}}R^{(b)})^{2p+1}}\, T_p(b,R^{(b)},\partial_{x_{b}}R^{(b)},\cdots,\partial^{p+1}_{x_{b}}R^{(b)}).
	\end{equation}
\end{lemma}
\begin{proof}
	Since $R\mapsto Z^{(b)}(R;x_{b+1},x_{b+2},\ldots)$ and $x_{b}\mapsto R^{(b)}(x_b,x_{b+1},\ldots)$ are compositional inverses (see \eqref{eq:Zfun}), it follows from Lagrange inversion (or more explicitly \cite[Theorem 1]{Johnson2002}) that
	\begin{equation*}
	(\partial_{x_b} R^{(b)})^{2k-1}\frac{\partial^k Z^{(b)}}{\partial R^k}( R^{(b)};x_{b+1},x_{b+2},\ldots)
	\end{equation*}
	is a homogeneous polynomial in $\partial_{x_{b}}R^{(b)},\cdots,\partial^{k}_{x_{b}}R^{(b)}$ of degree $k-1$.
	
	The expression for $\bar{M}_p^{(b)}$ in Proposition \ref{thm:highergenuspartitionfunctions} can also be written as
	\begin{equation*}
	\bar{M}_p^{(b)}(x_b,x_{b+1},\ldots) = r^{-b}\, Q_p(b,-b+r\, \partial_r)\, Z^{(b)}(r; x_{b+1}, x_{b+2}, \ldots) \Big|_{r=R^{(b)}(x_b,x_{b+1},\ldots)}.
	\end{equation*}
	Since $Q_p(-b,b)=0$ and $Q_p(b,j)$ is of degree $p+1$ in $j$ by Lemma \ref{lem:Qsums}, it corresponds to a linear expression in
	\begin{equation*}
	(R^{(b)})^{k-b}\frac{\partial^k Z^{(b)}}{\partial R^k}( R^{(b)};x_{b+1},x_{b+2},\ldots) \qquad \text{for }k=1,\ldots,p+1
	\end{equation*}
	with coefficients that are polynomials in $b$.
	Combining both statements, we conclude that
	\begin{equation*}
	(\partial_{x_b} R^{(b)})^{2p+1} (R^{(b)})^{b-1} \bar{M}_p^{(b)}(x_b,x_{b+1},\ldots)
	\end{equation*}
	is a homogeneous polynomial of degree $2p$ in  $R^{(b)}$ and its first $p+1$ derivatives with coefficients that are polynomials in $b$.
\end{proof}

The first few polynomials read
\begin{align*}
T_0 &= 1,\\
T_1 &= \frac{2}{3} b(b-1) r_1^2 - \frac{2}{3} r_0 r_2,\\
T_2 &= \frac{1}{30} (2b+1)b(b-1)(2b-3) r_1^4 - \frac{1}{30} (8b^2-16b+5)r_0 r_1^2r_2+\frac{4}{5} r_0^2 r_2^2 - \frac{4}{15} r_0^2r_1r_3,\\
T_3 &= \frac{1}{315} (b-2) (b-1) b (b+1) (2 b-3) (2 b+1) r_1^6-\frac{2}{105} (b-2) b (2 b^2-2 b+1) r_0 r_2 r_1^4\\
&\quad+
\frac{2}{35} (4 b^2-12 b+7) r_0^2r_1^2 r_2^2-\frac{2}{105} (4 b^2-12 b+7) r_0^2r_1^3 r_3-\frac{8
	}{7}r_2^3+\frac{16}{21}r_0^3 r_1 r_3 r_2-\frac{8}{105} r_0^3r_1^2 r_4.
\end{align*}
In particular, Proposition \ref{thm:highergenuspartitionfunctions} leads to the following explicit expressions for the generating functions of genus-$1$ and genus-$2$ essentially $2b$-irreducible maps with $b\geq 1$,
\begin{align*}
F_1^{(b)}(x_b,x_{b+1},\ldots) &= -\frac{1}{12}\log \frac{(R^{(b)})^{1-b}}{\partial_{x_b}R^{(b)}},\\
F_2^{(b)}(x_b,x_{b+1},\ldots) &= \frac{1}{240} +r_0^{2b}\Bigg[\frac{\left(34 b^2-32 b+3\right) r_2^2}{2880 r_1^2}+\frac{\left(43 b^4-62 b^3+50 b^2-19 b+6\right) r_2}{1440
	r_0}\\
&\quad+\frac{\left(-86 b^6+258 b^5-413 b^4+396 b^3-263 b^2+108 b-36\right) r_1^2}{8640 r_0^2}\\
&\quad-\frac{(3 b-1) (4 b-1) r_3}{720 r_1}+\frac{r_2^3r_0}{90 r_1^4}-\frac{7 r_3 r_2r_0}{480 r_1^3}+\frac{r_4r_0}{288
	r_1^2} \Bigg],\quad\quad r_k \coloneqq \frac{\partial^k R^{(b)}}{\partial x_b^k}.
\end{align*}

\section{Generating functions without vertices of degree one}\label{sec:hatgenfun}

Let $\hat{\mathcal{M}}^{(b)}_{g}\subset \mathcal{M}^{(b)}_{g}$ be the set of essentially $2b$-irreducible maps of genus-$g$ with no vertices of degree $1$ and $\hat{\vec{\mathcal{M}}}^{(b)}_{g}\subset \vec{\mathcal{M}}^{(b)}_{g}$ the set of such maps that are rooted by distinguishing an oriented edge.
For $b\geq 1$ the \emph{$2b$-irreducible genus-$g$ partition function} is the formal generating function
\begin{align*}
\hat{F}_g^{(b)}(\hat{x}_{b},\hat{x}_{b+1}, \ldots) &= \sum_{\map\in\hat{\mathcal{M}}^{(b)}_g} \frac{1}{|\operatorname{Aut}(\map)|} \prod_{f\in \mathcal{F}(\map)} \hat{x}_{\deg(f)/2}  = \sum_{\map\in\hat{\vec{\mathcal{M}}}^{(b)}_g} \frac{1}{2|\mathcal{E}(\map)|} \prod_{f\in \mathcal{F}(\map)} \hat{x}_{\deg(f)/2},
\end{align*}
where $\mathcal{F}(\map)$ is the set of faces of $\map$ and $|\mathcal{E}(\map)|$ is the number of edges.
We will also consider the case $b=0$ where $\hat{x}_0$ appears in the weight assigned to the vertices of the map and the partition function includes all genus-$g$ maps (with vertices of degree one allowed),
\begin{align*}
\hat{F}_g^{(0)}(\hat{x}_{0},\hat{x}_{2}, \ldots) &= \sum_{\map\in\mathcal{M}^{(0)}_g} \frac{1}{|\operatorname{Aut}(\map)|} \prod_{v\in \mathcal{V}(\map)} ( \hat{x}_0 + \ind_{\{\deg(v)\geq 2\}}) \prod_{f\in \mathcal{F}(\map)} \hat{x}_{\deg(f)/2}.
\end{align*}
Note that setting $\hat{x}_0=0$ gives the conventional generating function $\hat{F}_g(\hat{x}_1,\hat{x}_2,\ldots) = \hat{F}_g^{(0)}(0,\hat{x}_{0},\hat{x}_{1}, \ldots)$ of genus-$g$ maps without vertices of degree one.
It is useful to think of $F_g^{(0)}(x_0,x_1,\ldots)$ as the generating functions of maps in which the vertices are bi-colored, say in black and white, such that a black vertex carries weight $1$ and a white vertex weight $x_0$. 
The generating function $\hat{F}_g^{(0)}(\hat{x}_0,\hat{x}_1,\ldots)$ then has the same interpretation, except black vertices are required to have degree at least two.

With these definitions we have for any $b,g\geq 0$, $n\geq 1$ (provided $n\geq 3$ if $g=0$) and $\ell_1,\ldots,\ell_n \geq \max(b,1)$,
\begin{equation}\label{eq:numhatfrompartition}
\|\hat{\mathcal{M}}^{(b)}_{g,n}(\ell_1,\ldots,\ell_n)\| = \frac{\partial^n \hat{F}_g^{(b)}}{\partial \hat{x}_{\ell_1}\cdots\partial \hat{x}_{\ell_n}}(0,0,\ldots).
\end{equation}

\begin{proposition}\label{thm:hatsubstitution}
The partition functions $F^{(b)}_g$ and $\hat{F}^{(b)}_g$ for $g\geq 0$ and $b\geq 0$ are related via
\begin{align*}
F^{(b)}_g(x_b,x_{b+1},\ldots) &= \hat{F}^{(b)}_g(\hat{x}_b,\hat{x}_{b+1},\ldots)\Big|_{\hat{x}_p = \sum_{\ell=p}^\infty x_\ell A_{\ell,p}^{(b)}},\quad 
A^{(b)}_{\ell,p} = \ind_{\ell=p=b} + \frac{p}{\ell} \binom{2\ell}{\ell-p}\,\ind_{\ell \geq p > b},\\
\hat{F}^{(b)}_g(\hat{x}_b,\hat{x}_{b+1},\ldots) &= F^{(b)}_g(x_b,x_{b+1},\ldots)\Big|_{x_\ell = \sum_{p=\ell}^\infty \hat{x}_p B_{p,\ell}^{(b)}},\quad B_{p,\ell}^{(b)} = \ind_{p=\ell=b} + (-1)^{p-\ell} \binom{p+\ell-1}{p-\ell} \ind_{p\geq \ell> b}.
\end{align*}
\end{proposition}
\begin{figure}
	\centering
	\includegraphics[width=\linewidth]{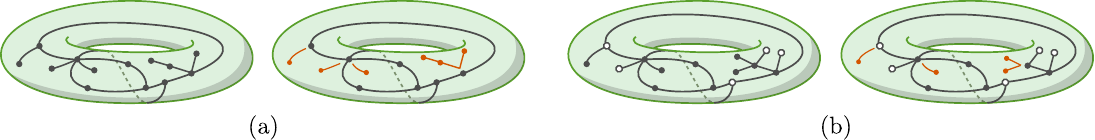}
	\caption{(a) A genus-$1$ map (left) and its core (right) together with the trees (in orange) associated to the corners of the core. (b) Similar but for a bi-colored genus-$1$ map.\label{fig:trees}}
\end{figure}
\begin{proof}
	We focus first on the case $b\geq 1$.
	The \emph{core} $\mathsf{Core}(\map)$ of a map $\map$ that is not a tree is the map obtained by repeatedly removing vertices of degree one and the incident edges (Figure \ref{fig:trees}a). 
	Given a rooted $2b$-irreducible genus-$g$ map $\map\in \hat{\vec{\mathcal{M}}}_g$ one may construct another map $\map'$ by inserting any number of trees into the corners of the map except for the corners that belong to a face of degree $2b$.
	Then $\map'$ will also be $2b$-irreducible since no new cycles have been produced and $\map = \mathsf{Core}(\map')$. 
	Vice versa, any map $\map'\in \vec{\mathcal{M}}_g$ that is rooted on an edge contained in its core $\mathsf{Core}(\map')$ can be obtained in this way from the map $\mathsf{Core}(\map')\in\hat{\vec{\mathcal{M}}}_g$.
	The number of ways one can insert trees into the corners of a face of degree $2p$ to obtain a face of degree $2\ell$ is
	\begin{equation}
	A^{(b)}_{\ell,p} = \ind_{\ell=p=b} + \frac{p}{\ell} \binom{2\ell}{\ell-p}\,\ind_{\ell \geq p > b}.
	\end{equation}
	Therefore
	\begin{align*}
	F_g^{(b)}(x_{b},x_{b+1}, \ldots) &= \sum_{\map'\in\vec{\mathcal{M}}^{(b)}_g} \frac{1}{2|\mathcal{E}(\map')|} \prod_{f\in \mathcal{F}(\map')} x_{\deg(f)/2}\\
	&= \sum_{\map'\in\vec{\mathcal{M}}^{(b)}_g} \frac{\ind_{\{\text{root in }\mathsf{Core}(\map')\}}}{2|\mathcal{E}(\mathsf{Core}(\map'))|} \prod_{f\in \mathcal{F}(\map')} x_{\deg(f)/2} \\
	&= \sum_{\map\in\hat{\vec{\mathcal{M}}}^{(b)}_g} \frac{1}{2|\mathcal{E}(\map)|} \prod_{f\in \mathcal{F}(\map)}\, \sum_{\ell = \deg(f)/2}^\infty A_{\ell,\deg(f)/2}^{(b)}\, x_\ell\\
	&= \hat{F}^{(b)}_g(\hat{x}_b,\hat{x}_{b+1},\ldots)\Big|_{\hat{x}_p = \sum_{\ell=b}^p x_\ell A_{\ell,p}^{(b)}}.
	\end{align*}
	
	In the case $b=0$, we make use of the interpretation of $F_g^{(0)}$ and $\hat{F}_g^{(0)}$ as generating functions of maps in which the vertices are bi-colored.
	We let the core $\mathsf{Core}(\map)$ of such a map $\map$ be defined as before, except that we do not remove white vertices of degree one (Figure \ref{fig:trees}b).
	Then the argument above works in this case as well: any map $\map'$ with bi-colored vertices is uniquely obtained from the map $\map = \mathsf{Core}(\map')$ without white vertices of degree one by inserting trees in the corners of $\map$.
	The same computation therefore verifies that
	\begin{equation*}
	F_g^{(0)}(x_{0},x_{1}, \ldots) = \hat{F}^{(0)}_g(\hat{x}_0,\hat{x}_{1},\ldots)\Big|_{\hat{x}_p = \sum_{\ell=0}^p x_\ell A_{\ell,p}^{(0)}}.
	\end{equation*}
	
	For the reverse substitution we note that the hypergeometric identity
	\begin{equation*}
	\sum_{p=k}^\ell \frac{p}{\ell}\binom{2\ell}{\ell-p}\, (-1)^{p-k} \binom{p+k-1}{p-k}\, = \ind_{\ell=k}
	\end{equation*}
	implies that for any $m\geq b\geq 0$, the matrix $(A_{\ell,p}^{(b)})_{\ell,p=b}^m$ is invertible and its inverse is given by $(B^{(b)}_{p,\ell})_{\ell,p=b}^m$.
\end{proof}

Since the partition functions $F_g^{(b)}$ given in Propositions \ref{thm:genus0partitionfunction} and \ref{thm:highergenuspartitionfunctions} are expressed almost entirely in terms of $R^{(b)}(x_b,x_{b+1},\ldots)$, we shall focus first on the analogue of $R^{(b)}$ when vertices of degree $1$ are suppressed.

\subsection{The building block $\hat{R}^{(b)}$}

It follows from Proposition \ref{thm:hatsubstitution} and \eqref{eq:Rinterpretation} that
\begin{equation}\label{eq:RhatvsR}
\hat{R}^{(b)}(\hat{x}_{b},\hat{x}_{b+1}, \ldots) \coloneqq R^{(b)}-1\Big|_{x_\ell = \sum_{p=\ell}^\infty \hat{x}_p B_{p,\ell}^{(b)}} =\frac{\partial^2F_0^{(b)}}{\partial x_{b+1}\partial x_{b}}\Big|_{x_\ell = \sum_{p=\ell}^\infty \hat{x}_p B_{p,\ell}^{(b)}} =  \frac{\partial^2\hat{F}_0^{(b)}}{\partial\hat{x}_{b+1}\partial\hat{x}_{b}}.
\end{equation}
Hence, if $b\geq 1$, $\hat{R}^{(b)}$ enumerates $2b$-irreducible planar maps with no vertices of degree one and two distinguished faces of degree $2b+2$ respectively $2b$.
In the case $b=0$ it is the generating function of planar maps with a distinguished face of degree $2$ and a marked vertex in which each non-marked vertex $v$ carries weight $\hat{x}_0 + \ind_{\{\deg(v)\geq 2\}}$.
To describe the analogue of equation \eqref{eq:Zfun} satisfied by $\hat{R}^{(b)}$ we introduce two formal power series that play a central role in the following.

Let $I(b,\ell;r)$ and $J(b;r)$ be the (hypergeometric) formal power series in $r$  defined by
\begin{align}
	I(b,\ell;r) &= \sum_{p=0}^\infty \frac{r^p}{(p!)^2}\prod_{m=0}^{p-1}(\ell^2-(b-m)^2) = {_2F_1}(\ell-b,-\ell-b;1;-r), \label{eq:Iseries}\\
	J(b;r) &= \sum_{p=1}^\infty \frac{(-1)^{p+1} r^p}{p!(p-1)!} \prod_{m=0}^{p-2} (b-m)(b-m-1) = r\,\,{_2F_1}(1-b,-b;2;-r),\label{eq:Jseries}
\end{align}
where ${_2F_1}$ is the hypergeometric function ${_2F_1}(a,b;c;z) = \sum_{n=0}^\infty \frac{(a)_n(b)_n}{(c)_n}\frac{z^n}{n!}$ in terms of the rising Pochhammer symbol $(a)_n = a(a+1)\cdots(a+n-1)$.
Clearly the coefficients of $I(b,\ell;r)$ are polynomials in $b$ and $\ell^2$ and those of $J(b;r)$ are polynomials in $b$.
The notations $I$ and $J$ are chosen because of the resemblance of their expansion to those of certain Bessel functions $I_0$ and $J_1$. 
Indeed, selecting the homogeneous part in $b$ and $\ell$ of top degree in the coefficients of $I(b,\ell;r)$ gives
\begin{equation}
I^{\text{hom}}(b,\ell;r) = \sum_{p=0}^\infty \frac{r^p}{(p!)^2} (\ell^2-b^2)^p = I_0\left(2\sqrt{(\ell^2-b^2)r}\right), 
\end{equation} 
where $I_0$ is a modified Bessel function of the first kind.
Similarly, selecting the top degree monomials in the coefficient of $J(b;r)$ yields
\begin{equation}
J^{\text{hom}}(b;r) = \sum_{p=1}^\infty \frac{(-1)^{p+1} r^p}{p!(p-1)!} b^{2p-2} = \frac{\sqrt{r}}{b} J_1(2b\sqrt{r}),
\end{equation}
where $J_1$ is a Bessel function of the first kind.
We will not use these relations here, but they play an important role in the follow-up work \cite{Budd2020}.

\begin{proposition}\label{thm:Rhat}
For $b\geq 0$ the formal power series $\hat{R}^{(b)}$ is the unique solution to $\hat{Z}^{(b)}(\hat{R}^{(b)};\hat{x}_{b},\hat{x}_{b+1},\ldots) =0$, where
\begin{equation}
\hat{Z}^{(b)}(r;\hat{x}_{b},\hat{x}_{b+1},\ldots) \coloneqq J(b,r) - \sum_{\ell=b}^\infty I(b,\ell;r)\,\hat{x}_\ell.
\end{equation}
\end{proposition}

\begin{proof}
We claim that $\hat{Z}^{(b)}$ and $Z^{(b)}$ are related via	
\begin{equation}\label{eq:ZhatvsZ}
	\hat{Z}^{(b)}(r;\hat{x}_b,\hat{x}_{b+1},\ldots) = -\hat{x}_b + Z^{(b)}(1+r;x_{b+1},x_{b+2})\Big|_{x_\ell = \sum_{p=\ell}^\infty \hat{x}_p B_{p,\ell}^{(b)}}.
\end{equation}
Combined with \eqref{eq:RhatvsR} and the fact that $R^{(b)}$ is uniquely determined by $Z^{(b)}(R^{(b)}; x_{b+1}, \ldots ) = x_b$, this proves the proposition.

To verify the identity we evaluate the right-hand side of \eqref{eq:ZhatvsZ} using \eqref{eq:Zfun} to
\begin{align*}
&(1+r)\ind_{\{b=0\}}-\sum_{\ell=0}^b (-1)^{b-\ell}\binom{b+\ell}{2\ell}\operatorname{Cat}(\ell)(1+r)^{b-\ell} - \hat{x}_b - \sum_{j=b+1}^\infty \binom{2j-1}{j+b}(1+r)^{b+j}\sum_{p=j}^\infty \hat{x}_p B_{p,j}^{(b)}\\
&= \tilde{J}(b,r) - \sum_{p=b}^\infty \tilde{I}(b,p;r)\,\hat{x}_p,
\end{align*}
where 
\begin{align}
\tilde{J}(b;r)&=(1+r)\ind_{\{b=0\}}-\sum_{\ell=0}^b (-1)^{b-\ell}\binom{b+\ell}{2\ell}\operatorname{Cat}(\ell)(1+r)^{b-\ell},\\
\tilde{I}(b,p;r)&=\ind_{\{b=p\}}+\sum_{\ell=b+1}^p (-1)^{p-\ell} \binom{p+\ell-1}{p-\ell} \binom{2\ell-1}{\ell+b} (1+r)^{\ell+b}.\label{eq:Isumidentity}
\end{align}
It remains to show that $\tilde{J}(b,r)=J(b,r)$ and $\tilde{I}(b,p;r)=I(b,p;r)$ for $p \geq b\geq 0$, which we could do by appealing to known hypergeometric identities.
Instead we choose to provide self-contained proofs involving generating functions and Lagrange inversion, since the intermediate expressions will be of use later on.
We note that $\tilde{J}(0;r) = J(0;r) = r$, while for $b\geq 1$ we may use the well-known generating function $\frac{1-\sqrt{1-4t}}{2t} = \sum_{\ell=0}^\infty \operatorname{Cat}(\ell) t^\ell$ of the Catalan numbers to compute the generating function
\begin{align*}
	U(z)\equiv U(r;z) &= 
\sum_{b=1}^\infty \tilde{J}(b;r) z^b \\
&=1-\sum_{\ell=0}^\infty \sum_{b=\ell}^\infty z^b (-1)^{b-\ell}\binom{b+\ell}{2\ell}\operatorname{Cat}(\ell)(1+r)^{b-\ell}\\
 &= 1-\sum_{\ell=0}^\infty \operatorname{Cat}(\ell) z^\ell (1+z+rz)^{-2\ell-1} \\
&= 1-(1+z+rz)^{-1} \frac{1-\sqrt{1-4t}}{2t}\Big|_{t = \frac{z}{(1+z+rz)^2}}\\ 
&= 1-\frac{1+z+r z - \sqrt{(1+z+r z)^2-4z}}{2z}.
\end{align*}
Since $U(r;z) = rz+ O(z^2)$ is the power series solution to
\begin{equation}
\frac{U}{(1-U)(r+U)} = z,
\end{equation}
Lagrange inversion \cite{Gessel2016} easily gives (see also \cite[Equation (5.16)]{Bouttier2014})
\begin{align*}
\tilde{J}(b;r)&= \frac{1}{b} [u^{b-1}] (1-u)^b(r+u)^b 
= \sum_{p=1}^\infty r^p \frac{1}{b}\binom{b}{p}[u^{b-1}](1-u)^bu^{b-p}
=\sum_{p=1}^\infty (-1)^{p-1}r^p \frac{1}{b}\binom{b}{p}\binom{b}{p-1}
\end{align*}
for $b\geq 1$.
Expanding the binomials gives precisely $J(b,r)$ in \eqref{eq:Jseries}.

Since $\tilde{I}(b,b;r)=I(b,b;r)=1$, it remains to demonstrate $\tilde{I}(b,p;r)=I(b,p;r)$ for $p>b\geq 0$.
To this end we examine the generating function 
\begin{align*}
	\sum_{p=b+1}^\infty \tilde{I}(b,p;r) y^p = \sum_{\ell=b+1}^\infty \binom{2\ell-1}{\ell+b} (1+r)^{\ell+b} \frac{y^\ell}{(1+y)^{2\ell}}.
\end{align*}
With the help of the identity
\begin{align*}
\sum_{\ell=b+1}^\infty \binom{2\ell-1}{\ell+b} t^\ell = \frac{t^{-b}}{\sqrt{1-4t}} \left(\frac{1-\sqrt{1-4t}}{2}\right)^{2b+1},
\end{align*}
and the rational parametrization 
\begin{align*}
	t=\frac{(1+r) y}{(1+y)^2}, \qquad y = \frac{\tilde{U}(1-\tilde{U})}{1+r-\tilde{U}}\qquad\text{such that}\quad \frac{1-\sqrt{1-4t}}{2} = \frac{\tilde{U}}{1+y},
\end{align*}
this reduces to
\begin{align}
	\sum_{p=b+1}^\infty \tilde{I}(b,p;r) y^p = \frac{\tilde{U}^{2b+1}}{y^b}\frac{1}{1+y - 2\tilde{U}}  = \frac{\tilde{U}^{b+1}(1+r-\tilde{U})^{b}}{(1-\tilde{U})^{b+1}}\frac{1}{1-\frac{r\tilde{U}}{(1+r-\tilde{U})(1-\tilde{U})}}.\label{eq:Igenfun}
\end{align}
Writing $\mathcal{R}(u) =(1+r-u)/(1-u)$, the formal power series 
\begin{align*}
	\tilde{U}(r;y) = \frac{1+y-\sqrt{(1-y)^2-4r y}}{2} = (r+1)y + O(y^2)
\end{align*}
is the unique solutions to $\tilde{U} = y\,\mathcal{R}(\tilde{U})$.
Writing $\psi(u) = u^{b+1}(1+r-u)^b(1-u)^{-b-1}$, Lagrange inversion in the form \cite[Equation (2.1.4)]{Gessel2016} then gives
\begin{align*}
\tilde{I}(b,\ell;r) &= [y^\ell]\frac{\psi(\tilde{U})}{1-y\, \mathcal{R}'(\tilde{U})}\\
&= [u^\ell]\psi(u) \mathcal{R}(u)^\ell\\
&= [u^{\ell-b-1}](1+r-u)^{b+\ell}(1-u)^{-b-\ell-1}\\
&= \sum_{p=0}^{\ell+b} r^p \binom{b+\ell}{p}\binom{p+\ell-b-1}{p}.
\end{align*}
Again, after expanding the binomials we recover $I(b,\ell;r)$ in \eqref{eq:Iseries}.
\end{proof}

\begin{lemma}\label{lem:Rder}
 For $\ell \geq b \geq 0$ we have the identity
	\begin{equation}
	\frac{\partial\hat{R}^{(b)}}{\partial \hat{x}_{\ell}}  = I(b,\ell;\hat{R}^{(b)}) \frac{\partial\hat{R}^{(b)}}{\partial \hat{x}_{b}}.
	\end{equation}
\end{lemma}
\begin{proof}
	This follows directly from taking derivatives of the equation $\hat{Z}^{(b)}(\hat{R}^{(b)};\hat{x}_{b},\hat{x}_{b+1},\ldots)=0$.
\end{proof}

The following result gives a compact expression for the coefficients of $\hat{R}^{(b)}$, allowing in particular to deduce that they are polynomials in $b$ and the square face degrees.
The result is analogous to formulas in \cite[Section 2.3-2.4]{Bouttier2014a}, that were obtained via combinatorial means, but here we give an algebraic proof.
To state the result we introduce the formal power series inverse 
\begin{equation*}
z \mapsto J^{-1}(b;z) = z + \tfrac{1}{2}b(b-1) z^2 + \tfrac{1}{12}b(b-1)^2(5b+2) z^3 + \cdots
\end{equation*}
to $r\mapsto J(b;r)$, whose coefficient of $z^k$ is easily seen to be a polynomial in $b$ of degree $2k-2$.

\begin{lemma}\label{lem:Rderpolynomial}
For $b\geq 0$, $k\geq 1$ and $\ell_1,\ldots,\ell_k \geq b$ and $f(r)$ any formal power series with $f'(0)\neq 0$,
\begin{equation}\label{eq:Rderexpr0}
\frac{\partial^k f(\hat{R}^{(b)})}{\partial \hat{x}_{\ell_1}\cdots \partial \hat{x}_{\ell_k}}(0,0,\ldots) = k!\, [z^k] \int_0^{J^{-1}(b;z)} \!\!\!\!\rmd r\,f'(r)\, \prod_{i=1}^k I(b,\ell_i;r)
\end{equation}
is a polynomial in $b, \ell_1^2, \ldots, \ell_k^2$ that is symmetric and of degree $k-1$ in $\ell_1^2, \ldots, \ell_k^2$.
\end{lemma}
\begin{proof}
We claim that for any $k\geq 0$,
\begin{equation}\label{eq:Rderexpr}
\frac{\partial^k f(\hat{R}^{(b)})}{\partial \hat{x}_{\ell_1}\cdots \partial \hat{x}_{\ell_k}} = f(0)\,\ind_{\{k=0\}}+\frac{\partial^k}{\partial\hat{x}_{b}^k} \int_0^{\hat{R}^{(b)}}\!\!\!\!\rmd r\, f'(r)\, \prod_{i=1}^k I(b,\ell_i;r),
\end{equation}
from which the required identity follows by evaluation at ${\hat{x}_{b} = \hat{x}_{b+1} = \cdots =0}$, where we use that $\hat{R}^{(b)} = J^{-1}(b,\hat{x}_{b})$ when $\hat{x}_{\ell}=0$ for all $\ell> b$.
We verify the claim by induction on $k$, noting that the base case $k=0$ is trivially satisfied.
If \eqref{eq:Rderexpr} holds for $k \geq 0$, then
\begin{align*}
\frac{\partial^{k+1} f(\hat{R}^{(b)})}{\partial \hat{x}_{\ell_1}\cdots \partial \hat{x}_{\ell_{k+1}}} &=\frac{\partial^k}{\partial \hat{x}_{b}^k}\,\frac{\partial}{\partial\hat{x}_{\ell_{k+1}}}\int_0^{\hat{R}^{(b)}}\!\!\!\!\rmd r\,f'(r)\, \prod_{i=1}^k I(b,\ell_i;r) 
= \frac{\partial^k}{\partial \hat{x}_{b}^k}\,\frac{\partial f(\hat{R}^{(b)})}{\partial\hat{x}_{\ell_{k+1}}} \prod_{i=1}^k I(b,\ell_i;\hat{R}^{(b)})\\
&\leftstackrel{\text{Lem. }\ref{lem:Rder}}{=}\frac{\partial^k}{\partial \hat{x}_{b}^k}\frac{\partial f(\hat{R}^{(b)})}{\partial\hat{x}_{b}} \prod_{i=1}^{k+1} I(b,\ell_i;\hat{R}^{(b)})
= \frac{\partial^{k+1}}{\partial \hat{x}_{b}^{k+1}}\int_0^{\hat{R}^{(b)}}\!\!\!\!\rmd r\,f'(r)\, \prod_{i=1}^{k+1} I(b,\ell_i;r),
\end{align*}
verifying the induction step.

The power series coefficients of $J^{-1}(b;z)$ and $I(b,\ell_i;r)$ are polynomials in $b,\ell_1^2,\ldots,\ell_k^2$ and therefore the same is true for the right-hand side of \eqref{eq:Rderexpr0}.
To see that it is of degree $k-1$ in $\ell_1^2,\ldots,\ell_k^2$, we note that the top-degree monomials arise from the contribution
\begin{equation*}
k! [z^k] \int_0^z \rmd r\,f'(0) \prod_{i=1}^k I(b,\ell_i;r) = f'(0) (k-1)! [r^{k-1}] \prod_{i=1}^k I(b,\ell_i;r)
\end{equation*}
and use that the power series coefficient $[r^p]I(b,\ell_i;r)$ is polynomial $\ell_i^2$ of degree $p$.
\end{proof}

\subsection{Generating functions for arbitrary genus}

All ingredients are in place to perform the substitution of Proposition~\ref{thm:hatsubstitution} in the partition functions of Proposition~\ref{thm:genus0partitionfunction} and Proposition~\ref{thm:hatsubstitution}.

\begin{proposition}\label{thm:parthat}
For $b\geq 0$ the partition functions of essentially $2b$-irreducible maps with no vertices of degree one are given by (with $\ell_1,\ell_2\geq \max(b,1)$)
\begin{align}
\frac{\partial}{\partial \hat{x}_{\ell_1}}\frac{\partial}{\partial \hat{x}_{\ell_2}}\hat{F}_0^{(b)}(\hat{x}_b,\hat{x}_{b+1},\ldots) &= \frac{3\hat{x}_b^2+2\hat{x}_b^3}{2(1+\hat{x}_b)^2}\,\ind_{\ell_1=\ell_2=b} \, +  \int_{\hat{x}_0\ind_{\{b=0\}}}^{\hat{R}^{(b)}} \!\!\!\rmd r\, \frac{I(b,\ell_1;r)I(b,\ell_2;r)}{(1+r)^{2b+1}}, \label{eq:Fhat0part}\\
\hat{F}^{(b)}_1(\hat{x}_b,\hat{x}_{b+1},\ldots) &= - \frac{1}{12} \log \hat{M}^{(b)}_0\mkern60mu\text{for }b\geq 1, \label{eq:Fhat1part}\\
\hat{F}^{(b)}_g(\hat{x}_b,\hat{x}_{b+1},\ldots) &= P_g\left( \frac{1}{\hat{M}^{(b)}_0}, \frac{\hat{M}^{(b)}_1}{\hat{M}^{(b)}_0},\ldots, \frac{\hat{M}^{(b)}_{3g-3}}{\hat{M}^{(b)}_0}\right) \quad \text{for }b\geq 1,\,g\geq 2,\label{eq:Fhatgpart}\\
\hat{F}^{(0)}_1(\hat{x}_0,\hat{x}_{1},\ldots) &= - \frac{1}{12} \log\frac{ \hat{M}^{(0)}_0}{1+\hat{x}_0},\nonumber\\
\hat{F}^{(0)}_g(\hat{x}_0,\hat{x}_{1},\ldots) &= (1+\hat{x}_0)^{2-2g}P_g\left( \frac{1}{\hat{M}^{(0)}_0}, \frac{\hat{M}^{(0)}_1}{\hat{M}^{(0)}_0},\ldots, \frac{\hat{M}^{(0)}_{3g-3}}{\hat{M}^{(0)}_0}\right) \quad \text{for }g\geq 2,\nonumber
\end{align}
where
\begin{equation}\label{eq:mhat}
\hat{M}_p^{(b)} = Q_p(b,(1+r)\, \partial_r)\,(1+r)^{-b} \, \hat{Z}^{(b)}(r; \hat{x}_{b},\hat{x}_{b+1}, \ldots) \Big|_{r=\hat{R}^{(b)}(\hat{x}_b,\hat{x}_{b+1},\ldots)}.
\end{equation}
\end{proposition}
\begin{proof}
We start with the genus-$0$ case. Proposition \ref{thm:genus0partitionfunction} and Proposition \ref{thm:hatsubstitution} together imply that
\begin{align*}
\frac{\partial^2\hat{F}_0^{(b)}}{\partial \hat{x}_{p_1}\partial \hat{x}_{p_2}} &= \sum_{\ell_1=p_1}^\infty \sum_{\ell_2=p_2}^\infty B^{(b)}_{p_1,\ell_1}B^{(b)}_{p_2,\ell_2} \frac{\partial^2F_0^{(b)}}{\partial x_{\ell_1}\partial x_{\ell_2}}\\
&=\frac{3\hat{x}_b^2+2\hat{x}_b^3}{2(1+\hat{x}_b)^2}\,\ind_{\ell_1=\ell_2=b} \, +  \int_{1+\hat{x}_0\ind_{\{b=0\}}}^{\hat{R}^{(b)}+1} \!\!\!\rmd r\, \frac{\prod_{i=1}^2\left[\sum_{\ell_i=p_i}^\infty B^{(b)}_{p_i,\ell_i}\left(\ind_{\{\ell_i=b\}} + \binom{2\ell_i-1}{\ell_i+b} r^{\ell_i+b}\right)\right]}{r^{2b+1}}.
\end{align*}
Using \eqref{eq:Isumidentity} and shifting the integration variable by one, reproduces \eqref{eq:Fhat0part}.

In the light of Proposition \ref{thm:highergenuspartitionfunctions}, to verify the formulas for the higher-genus partition functions it suffices to check that 
\begin{equation*}
\hat{M}_p^{(b)}(\hat{x}_b,\hat{x}_{b+1},\ldots) = \bar{M}_p^{(b)}(x_b,x_{b+1},\ldots)\Big|_{x_\ell = \sum_{p=\ell}^\infty \hat{x}_p B_{p,\ell}^{(b)}}
\end{equation*}
is given correctly by \eqref{eq:mhat}, which follows easily from \eqref{eq:ZhatvsZ} and $Q_{p}(b,-b)=0$ (Lemma \ref{lem:Qsums}).
\end{proof}

In the planar case we can use this to find an expression for the coefficients of $\hat{F}_0^{(b)}$ that treats all faces on equal footing and is the analogue of \cite[Equation (3.2)]{Bouttier2014a} when vertices of degree one are suppressed.

\begin{proposition}\label{thm:numplanar}
For any $b\geq 0$ the number of planar $2b$-irreducible maps with no vertices of degree one and $n\geq 3$ labeled faces of degrees $2\ell_1,\ldots,2\ell_n \geq \max(2b,2)$ satisfies
\begin{align}
\|\hat{\mathcal{M}}^{(b)}_{0,n}(\ell_1,\ldots,\ell_n)\| &= \hat{N}_{0,n}^{(b)}(\ell_1, \ldots, \ell_n) \,+\, \ind_{\{n\geq 4,\,\ell_1=\cdots=\ell_n=b\}}\,  \frac{(n-1)!}{2}(-1)^n, \label{eq:M0nenumeration}\\
\hat{N}_{0,n}^{(b)}(\ell_1,\ldots,\ell_n) &= (n-2)! [z^{n-2}] \int_0^{J^{-1}(b;z)}\!\!\rmd r\, \frac{\prod_{i=1}^n I(b,\ell_i;r)}{(1+r)^{2b+1}}.\label{eq:Nhat0nexpr}
\end{align}
In particular, $\hat{N}_{0,n}^{(b)}(\ell_1,\ldots,\ell_n)$ is a polynomial in $b$, $\ell_1^2, \ldots,\ell_n^2$ that is symmetric in $\ell_1^2, \ldots,\ell_n^2$ of degree $n-3$.
\end{proposition}
\begin{proof}
Recall from \eqref{eq:numhatfrompartition} that
\begin{equation*}
\|\hat{\mathcal{M}}^{(b)}_{0,n}(\ell_1,\ldots,\ell_n)\| = \frac{\partial}{\partial \hat{x}_{\ell_1}}\cdots \frac{\partial}{\partial \hat{x}_{\ell_n}} \hat{F}_0^{(b)} \Big|_{\hat{x}_{b}=\hat{x}_{b+1}=\cdots=0},
\end{equation*}
whose value we should be able to deduce from \eqref{eq:Fhat0part}.
A simple calculation shows that for $n\geq 3$,
\begin{equation*}
\left(\frac{\partial}{\partial \hat{x}_b}\right)^{n-2} \frac{3\hat{x}_b^2+2\hat{x}_b^3}{2(1+\hat{x}_b)^2}\,\Bigg|_{\hat{x}_{b}=0} = (-1)^n\frac{(n-1)!}{2}\ind_{\{n\geq 4\}},
\end{equation*}
so \eqref{eq:M0nenumeration} is satisfied with
\begin{equation*}
\hat{N}_{0,n}^{(b)}(\ell_1,\ldots,\ell_n) = \frac{\partial}{\partial \hat{x}_{\ell_3}}\cdots \frac{\partial}{\partial \hat{x}_{\ell_n}} \int_0^{\hat{R}^{(b)}} \!\!\!\rmd r\, \frac{I(b,\ell_1;r)I(b,\ell_2;r)}{(1+r)^{2b+1}}\,\Bigg|_{\hat{x}_{b}=\hat{x}_{b+1}=\cdots=0}.
\end{equation*}
Applying Lemma \ref{lem:Rderpolynomial} to the power series
\begin{equation}
f(r) = \int_0^r \rmd r \frac{I(b,\ell_1;r)I(b,\ell_2;r)}{(1+r)^{2b+1}} = r + O(r^2),
\end{equation}
we conclude that $\hat{N}_{0,n}^{(b)}$ is a polynomial in $b,\ell_1^2,\ldots,\ell_n^2$ that is symmetric in $\ell_1^2,\ldots,\ell_n^2$ of degree $n-3$ and can be expressed as \eqref{eq:Nhat0nexpr}.
\end{proof}

For the higher-genus case, it is convenient to make use of the expressions for the moments derived in Lemma \ref{lem:momentRder}. 
After substitution these become
\begin{equation}\label{eq:MhatRexpr}
\hat{M}^{(b)}_p(\hat{x}_{b},\hat{x}_{b+1},\ldots) = \frac{(1+\hat{R}^{(b)})^{1-b}}{(\partial_{\hat{x}_{b}}\hat{R}^{(b)})^{2p+1}}\, T_p(b,1+\hat{R}^{(b)},\partial_{\hat{x}_{b}}\hat{R}^{(b)},\cdots,\partial^{p+1}_{\hat{x}_{b}}\hat{R}^{(b)}).
\end{equation}

\begin{proposition}\label{thm:polynomialhighergenus}
	For any $g,n\geq 1$, $b\geq 0$ and $\ell_1,\ldots,\ell_n\geq \max(b,1)$ we have 
	\begin{equation*}
	\|\hat{\mathcal{M}}^{(b)}_{g,n}(\ell_1,\ldots,\ell_n)\| = \hat{N}_{g,n}^{(b)}(\ell_1, \ldots, \ell_n),
	\end{equation*}
	where 
	\begin{equation}
	\hat{N}^{(b)}_{g,n}(\ell_1,\ldots,\ell_n) = \frac{\partial^n \hat{F}^{(b)}_{g}}{\partial \hat{x}_{\ell_1}\cdots \partial \hat{x}_{\ell_n}}(0,0,\ldots)\label{eq:fhatcoeff}
	\end{equation}
	is a polynomial in $b,\ell_1^2, \ldots, \ell_n^2$ that is symmetric in $\ell_1^2, \ldots, \ell_n^2$.
\end{proposition}
\begin{proof}
From the definition of $\bar{M}_p^{(b)}$ in \eqref{eq:Mbarbdef} we deduce that $\bar{M}_p^{(b)}(0,0,\ldots) = \ind_{\{p=0\}}$ since $U^{(b)}_k(0,0,\ldots)=0$ for all $k\geq 1$.
Hence we also have $\hat{M}_p(0,\ldots) = T_p(b,1,0,\ldots,0) = \ind_{p=0}$.
Note further that $\hat{R}^{(b)}(0,\ldots) = 0$, $ \partial_{\hat{x}_{b}}\hat{R}^{(b)}(0,\ldots) = 1$.
It then follows easily from the expressions \eqref{eq:Fhat1part}, \eqref{eq:Fhatgpart} and \eqref{eq:MhatRexpr} that $\hat{N}^{(b)}_{g,n}(\ell_1,\ldots,\ell_n)$ is a polynomial in the variables 
\begin{equation*}
\left(\prod_{i\in I} \frac{\partial}{\partial\hat{x}_{\ell_i}}\right)\frac{\partial^k}{\partial\hat{x}_{b}^k} \hat{R}^{(b)}\Bigg|_{\hat{x}_b=\hat{x}_{b+1}=\ldots=0}, \qquad I\subset\{1,\ldots,n\}, \quad k\geq 0. 
\end{equation*}
Since each of these is polynomial in $\ell_1^2, \ldots, \ell_n^2$ by Lemma \ref{lem:Rderpolynomial}, the same is true for the power series coefficient $\hat{N}^{(b)}_{g,n}(\ell_1,\ldots,\ell_n)$ in \eqref{eq:fhatcoeff}.
The latter is symmetric in $\ell_1, \ldots, \ell_n$ by construction.
Together with \eqref{eq:numhatfrompartition} this concludes the proof.
\end{proof}

\subsection{Proof of Theorem \ref{thm:mainpolynomial}}

We have already established all ingredients for the proof of Theorem \ref{thm:mainpolynomial}. 
The first statement on the existence of a symmetric polynomial $\hat{N}_{g,n}^{(b)}$ satisfying \eqref{eq:mainenumeration} has been asserted in Proposition \ref{thm:numplanar} for $g=0$ and Proposition \ref{thm:polynomialhighergenus} for $g\geq 1$.
These propositions also demonstrate that the dependence of $\hat{N}_{g,n}^{(b)}$ on $b$ is polynomial.
It follows from Proposition \ref{thm:quasipolynomial} that the degree of the polynomial $\hat{N}_{g,n}^{(b)}$ both in $\ell_1,\ldots,\ell_n$ and in $b,\ell_1,\ldots,\ell_n$ is $2n+6g-6$.
Finally, the relation between $\|\mathcal{M}^{(b)}_{g,n}(\ell_1,\ldots,\ell_n)\|$ and $\|\hat{\mathcal{M}}^{(b)}_{g,n}(\ell_1,\ldots,\ell_n)\|$ in \eqref{eq:MvsMhat} is a direct consequence of \eqref{eq:numfrompartition}, \eqref{eq:numhatfrompartition} and Proposition \ref{thm:hatsubstitution}.\qed

\section{String and dilaton equations} \label{sec:stringdilaton}

The string and dilaton equations will be seen to be a consequence of the following special properties of the formal power series $I(b,\ell;r)$ and $J(b;r)$.

\begin{lemma}\label{thm:IJprop}
	The formal power series $I(b,\ell;r)$ and $J(b;r)$ satisfy
	\begin{enumerate}[label = {(\roman*)}]
		\item\label{item:IJrel} $I(b,1;r) = (1+r)^{b+1} \frac{\partial}{\partial r}((1+r)^{-b} J(b;r))$ for $b\geq 0$,
		\item\label{item:IJrel2} $I(b,1;r) - I(b,0;r) = J(b;r)$ for $b\geq 0$,
		\item\label{item:Isum} $(1+r)^{b+1} \tfrac{\partial}{\partial r} (1+r)^{-b}I(b,\ell;r) = - \ell\,I(b,\ell;r)+\sum_{k=1}^{\ell} 2k\, I(b,k;r) - \sum_{k=1}^{b} 2k\, I(b,k;r)$ for $\ell,b\geq 0$.
	\end{enumerate}
\end{lemma}
\begin{proof}
Recall from the proof of Proposition \ref{thm:Rhat} that 
\begin{align*}
I(b,\ell;r) &= [u^{b+\ell}] (1-u)^{b-\ell}(r+u)^{b+\ell}, \qquad (b+\ell\geq 0)\\
J(b;r) &= \frac{1}{b}[u^{b-1}] (1-u)^b(r+u)^b. \qquad\mkern22mu(b \geq 1)
\end{align*}
Note that we only proved the first relation for $\ell > b \geq 1$, but since the coefficients of $r$ on both sides are polynomial in $b$ and $\ell$ it must be valid for any $b,\ell\in\Z$ such that $b+\ell \geq 0$.

Property \ref{item:IJrel} is clearly satified in case $b=0$ since $I(0,1;r)=1+r$ and $J(0,r)=r$.
For $b\geq 1$ it follows from examining the coefficient of $u^{b-1}$ in
\begin{equation*}
(1+r)^{b+1} \frac{\partial}{\partial r}\left((1+r)^{-b} \frac{1}{b}(1-u)^b(r+u)^b\right) = (1-u)^{b+1}(r+u)^{b-1}
\end{equation*}
and the fact that $I(b,-\ell;r) = I(b,\ell;r)$.

For Property \ref{item:IJrel2} we use the identity 
\begin{equation*}
	\prod_{m=0}^{p-1} (1-b+m)(1+b-m) - \prod_{m=0}^{p-1} (-b+m)(b-m) = (-1)^{p-1} p \prod_{m=0}^{p-2} (b-m-1)(b-m),
\end{equation*}
which is easily checked by noting that the three terms have most factors in common.
The required identity then follows directly from the definitions \eqref{eq:Iseries} and \eqref{eq:Jseries} by comparing the coefficients.

The coefficients of $r$ on the both of the identity in Property \ref{item:Isum} are polynomials in $\ell$ and $b$, so it is sufficient to check the equality for $\ell > b \geq 0$.
To this end we rewrite \eqref{eq:Igenfun} as
\begin{align*}
K_b(r;y) &\coloneqq (1+r)^{-b}\sum_{p=b+1}^\infty I(b,p;r) y^p = \left(\frac{\tilde{U}(1+r-\tilde{U})}{(1+r)(1-\tilde{U})}\right)^b \frac{\tilde{U}}{1-\tilde{U}-\frac{r\tilde{U}}{1+r-\tilde{U}}}\\
&=\left(\frac{1-T}{1+T}\right)^b\frac{1-T}{2T}, \qquad T(r;y) = \frac{1}{1+y}\sqrt{(1-y)^2-4ry}. 
\end{align*}
We may easily check that
\begin{equation*}
(1+r)\frac{\partial T}{\partial r} = \frac{1+y}{1-y}\,y\, \frac{\partial T}{\partial y},
\end{equation*} 
implying that $K_b(r;y)$ satisfies the same first-order linear differential equation.
Hence
\begin{align*}
\sum_{p=b+1}^\infty (1+r)^{b+1} \frac{\partial}{\partial r}((1+r)^{-b} I(b,p;r))\, y^p &= (1+r)^{b+1}\frac{\partial K_b}{\partial r} =(1+r)^{b} \frac{1+y}{1-y}\,y\, \frac{\partial K_b}{\partial y}\\
&= \sum_{p=b+1}^\infty I(b,p;r) \frac{1+y}{1-y}\,p\, y^p.
\end{align*}
Extracting the coefficient of $y^\ell$, $\ell > b$, on both sides gives precisely the claimed identity.
\end{proof}

We will use these properties to prove the string and dilaton equations for essentially $2b$-irreducible genus-$g$ maps. 
We treat the planar case and the higher genus case separately, because of the differences in the expressions we have for their partition functions listed in Proposition \ref{thm:parthat}. 

\subsection{The planar case}

In the planar case we are stuck with having (at least) two distinguished faces in the generating function, so if we want to treat all faces on equal footing it is more convenient to use directly the expression \eqref{eq:Nhat0nexpr} for the polynomial $\hat{N}^{(b)}_{0,n}$.

\begin{proposition}\label{thm:stringdilatonplanar}
	The polynomials $\hat{N}_{0,n}^{(b)}(\ell_1,\ldots,\ell_n)$ satisfy the string and dilaton equations for $n\geq3$ and $b\geq 0$,
	\begin{align*}
		\hat{N}^{(b)}_{0,n+1}(\ell_1,\ldots,\ell_n,1) &= \sum_{j=1}^n \sum_{k=b+1}^{\ell_j} 2 k\, \hat{N}^{(b)}_{0,n}(\ell_1,\ldots,\ell_{j-1},k,\ell_{j+1},\ldots,\ell_n) - \sum_{j=1}^n \ell_j \hat{N}^{(b)}_{0,n}(\ell_1,\ldots,\ell_n),\\
		\hat{N}_{0,n+1}^{(b)}(\ell_1,\ldots,\ell_n,1) &- \hat{N}_{0,n+1}^{(b)}(\ell_1,\ldots,\ell_n,0) = (n-2) \hat{N}_{0,n}^{(b)}(\ell_1,\ldots,\ell_n).
	\end{align*}
\end{proposition}
\begin{proof}
	Using the shorthand notations
	\begin{equation}\label{eq:shorthand}
	[z^{n-1}] f(r) = [z^{n-1}]f(J^{-1}(b;z)), \quad \int \rmd r \,f(r) = \int_0^r \rmd s\,f(s), \quad J'(b;r) = \tfrac{\partial}{\partial r}J(b;r),
	\end{equation}
	we observe that any formal power series $f(r)$ satisfies
	\begin{align}
	(n-1)[z^{n-1}]\int\rmd r \,f(r) J'(b;r) &= (n-1)[z^{n-1}]\int\rmd z \,f(r) = [z^{n-2}]f(r) = [z^{n-2}]\int\rmd r\,f'(r),\\
	(n-1)[z^{n-1}]\int\rmd r \,f(r)J(b;r) &= (n-1)[z^{n-1}]\left(J(b;r) \int \rmd r \,f(r) -\int\rmd r\, J'(b;r) \int \rmd r \,f(r)\right)\nonumber\\ 
	&= (n-1)[z^{n-2}]\int \rmd r \,f(r) - [z^{n-2}]\int \rmd r \,f(r) \nonumber\\
	&= (n-2)[z^{n-2}]\int \rmd r \,f(r).\label{eq:Jinintegral0}
	\end{align}
	More generally we thus have for any formal power series $f(r)$ and $h(r)$ that
	\begin{equation}\label{eq:Jinintegral}
	(n-1)[z^{n-1}]\int\rmd r\,f(r) \tfrac{\partial}{\partial r} ( h(r) J(b;r) ) = [z^{n-2}]\int\rmd r \left( (n-1)f(r)h'(r) + f'(r) h(r) \right).
	\end{equation}
	According to Proposition \ref{thm:numplanar} we have for $n\geq3$ that
	\begin{align*}
	\hat{N}_{0,n+1}^{(b)}(1,\ell_1,\ldots,\ell_n) &\stackrel{\eqref{eq:Nhat0nexpr}}{=}(n-1)! [z^{n-1}]\int \rmd r\, (1+r)^{-1-2b}  I(b,1;r)\prod_{i=1}^nI(b,\ell_i;r),
	\end{align*}
	which with the help of Lemma~\ref{thm:IJprop}\ref{item:IJrel} and \eqref{eq:Jinintegral} gives
	\begin{align*}
	\hat{N}_{0,n+1}^{(b)}(1,\ell_1,\ldots,\ell_n)&= (n-1)! [z^{n-1}]\int \rmd r \, (1+r)^{-b} (\tfrac{\partial}{\partial r} (1+r)^{-b}J(b;r))\prod_{i=1}^nI(b,\ell_i;r)  \\
	&\stackrel{\eqref{eq:Jinintegral}}{=} (n-2)![z^{n-2}] \int \rmd r \Bigg( (n-1) (\tfrac{\partial}{\partial r} (1+r)^{-b}) (1+r)^{-b} \prod_{i=1}^n I(b,\ell_i;r)\\
	& \qquad\qquad\qquad\qquad\qquad + (1+r)^{-b} \tfrac{\partial}{\partial r} (1+r)^{-b} \prod_{i=1}^n I(b,\ell_i;r) \Bigg) \\
	&= (n-2)![z^{n-2}] \int \rmd r (1+r)^{-b} \sum_{j=1}^n \left(\prod_{i\neq j} I(b,\ell_i;r)\right)\tfrac{\partial}{\partial r} (1+r)^{-b} I(b,\ell_j;r).
	\end{align*}
	Applying Lemma~\ref{thm:IJprop}\ref{item:Isum} we thus find
	\begin{align*}
	\hat{N}_{0,n+1}^{(b)}(1,\ell_1,\ldots,\ell_n)&= (n-2)![z^{n-2}] \int \rmd r (1+r)^{-1-2b} \\
	&\qquad\qquad \times \sum_{j=1}^n \left(\prod_{i\neq j} I(b,\ell_i;r)\right)\left(-\ell_j I(b,\ell_j;r)+\sum_{k=b+1}^{\ell_j} 2k\, I(b,k;r)\right) \\
	&\stackrel{\eqref{eq:Nhat0nexpr}}{=} \sum_{j=1}^n \sum_{k=b+1}^{\ell_j} 2k\, \hat{N}_{0,n}^{(b)}(\ell_1,\ldots,\ell_{j-1},k,\ell_{j+1},\ldots,\ell_n) - \sum_{j=1}^n\ell_j \hat{N}_{0,n}^{(g)}(\ell_1,\ldots,\ell_n).	
	\end{align*}
	For the dilaton equation we use Lemma \ref{thm:IJprop}\ref{item:IJrel2} to calculate
	\begin{align*}
	\hat{N}_{0,n+1}^{(b)}(\ell_1,\ldots,&\ell_n,1) - \hat{N}_{0,n+1}^{(b)}(\ell_1,\ldots,\ell_n,0) \\
	&\stackrel{\eqref{eq:Nhat0nexpr}}{=} (n-1)! [z^{n-1}]\int \rmd r \left( I(b,1;r)-I(b,0;r) \right)(1+r)^{-1-2b} \prod_{i=1}^nI(b,\ell_i;r)\\
	&\stackrel{\mathmakebox[\widthof{=}]{\text{Lem.~\ref{thm:IJprop}\ref{item:IJrel2}}}}{=}\quad (n-1)! [z^{n-1}]\int \rmd r\,J(b;r) (1+r)^{-1-2b}\prod_{i=1}^nI(b,\ell_i;r)\\
	&\stackrel{\eqref{eq:Jinintegral0}}{=} (n-2)\, (n-2)! [z^{n-2}]\int \rmd r (1+r)^{-1-2b} \prod_{i=1}^nI(b,\ell_i;r) \\
	&\stackrel{\eqref{eq:Nhat0nexpr}}{=} (n-2)\hat{N}_{0,n}^{(b)}(\ell_1,\ldots,\ell_n),
	\end{align*}
	valid for any $n\geq 3$, which concludes the proof.
	\end{proof}

\subsection{The higher-genus case}

For genus $1$ and higher we do not have an explicit formula for the polynomials $\hat{N}^{(b)}_{g,n}$ like in the planar case.
The natural route to take is therefore to reformulate the string and dilation equations on the level of generating functions.
This however poses an immediate problem, because the equations involve an evaluation where one of the face (half-)degrees is set to $0$ or $1$, which does not make sense combinatorially in a $2b$-irreducible map when $b > 1$.

Luckily, we can salvage the situation as follows. 
Let $\check{R}^{(b)}(\hat{x}_0,\hat{x}_1,\ldots)$ be the power series solution to $\check{Z}^{(b)}(\check{R}^{(b)};\hat{x}_0,\hat{x}_1,\ldots)=0$ where 
\begin{equation}\label{eq:Zcheckdef}
\check{Z}^{(b)}(r; \hat{x}_{0}, \hat{x}_1, \ldots) = J(b, r) - \sum_{\ell=0}^{\infty} I(b,\ell;r) \hat{x}_{\ell},
\end{equation}
i.e.\ we have extended the definition of $\hat{Z}^{(b)}$ by including generating variables $\hat{x}_0, \ldots \hat{x}_{b-1}$.
Then we let $\check{F}_g^{(b)}(\hat{x}_0,\hat{x}_1,\ldots)$ for $g\geq 1$ and $\check{M}_p(\hat{x}_0,\hat{x}_1,\ldots)$ for $p\geq 0$ be defined by expressions \eqref{eq:Fhat1part}-\eqref{eq:mhat} of Proposition \ref{thm:parthat}, but with $\hat{R}^{(b)}$ and $\hat{Z}^{(b)}$ replaced by $\check{R}^{(b)}$ and $\check{Z}^{(b)}$.
Lemmas \ref{lem:Rder} and \ref{lem:Rderpolynomial} apply equally well to $\check{R}^{(b)}$, while Proposition \ref{thm:polynomialhighergenus} extends to the identity
\begin{equation}\label{eq:NfromFcheck}
	\hat{N}^{(b)}_{g,n}(\ell_1,\ldots,\ell_n) = \frac{\partial^n \check{F}^{(b)}_{g}}{\partial \hat{x}_{\ell_1}\cdots \partial \hat{x}_{\ell_n}}(0,0,\ldots),
\end{equation}
now valid for arbitrary non-negative $\ell_1,\ldots,\ell_n$.

\begin{lemma}\label{lem:strdilR}
Introducing the first-order linear partial differential operators
	\begin{align*}
	D^{\mathrm{str}} &= \frac{\partial}{\partial \hat{x}_1} - \sum_{\ell=0}^\infty \hat{x}_{\ell}\left( - \ell \frac{\partial}{\partial \hat{x}_{\ell}} +\sum_{k=1}^{\ell} 2k \frac{\partial}{\partial \hat{x}_{k}}-\sum_{k=1}^{b} 2k \frac{\partial}{\partial \hat{x}_{k}}\right), \\
	D^{\mathrm{dil}} &= \frac{\partial}{\partial \hat{x}_1} - \frac{\partial}{\partial \hat{x}_0} - \sum_{\ell=0}^\infty \hat{x}_{\ell} \frac{\partial}{\partial \hat{x}_{\ell}},
	\end{align*}
	we have
	\begin{equation}
	D^{\mathrm{str}} \check{R}^{(b)} = 1+\check{R}^{(b)}, \qquad D^{\mathrm{dil}} \check{R}^{(b)} = 0.
	\end{equation}
\end{lemma}
\begin{proof}
With the help of Lemmas \ref{lem:Rder} and \ref{thm:IJprop}, we may evaluate
	\begin{align*}
	D^{\mathrm{str}} \check{R}^{(b)} &= \left( I(b,1;\check{R}^{(b)}) - \sum_{\ell=0}^{\infty} \hat{x}_{\ell}\left[-\ell I(b,\ell;\check{R}^{(b)})  + \sum_{k=0}^\ell 2k I(b,k;\check{R}^{(b)}) - \sum_{k=0}^b 2k I(b,k;\check{R}^{(b)})\right]\right) \frac{\partial\check{R}^{(b)}}{\partial {\hat{x}_{b}}}\\
	&\mkern-30mu\stackrel{\text{Lem. \ref{thm:IJprop}\ref{item:IJrel},\ref{item:Isum}}}{=}  (1+\check{R}^{(b)})^{b+1}\frac{\partial\check{R}^{(b)}}{\partial {\hat{x}_{b}}} \frac{\partial}{\partial r} (1+r)^{-b} \left( J(b;r) - \sum_{\ell=0}^\infty \hat{x}_{\ell} I(b,\ell;r)\right)\Bigg|_{r=\check{R}^{(b)}} \\
	&=1+\check{R}^{(b)}+(1+\check{R}^{(b)})^{b+1}\frac{\partial}{\partial {\hat{x}_{b}}} (1+\check{R}^{(b)})^{-b} \left( J(b;\check{R}^{(b)}) - \sum_{\ell=0}^\infty \hat{x}_{\ell} I(b,\ell;\check{R}^{(b)})\right)\\
	&= 1+ \check{R}^{(b)},
	\end{align*}
	where in the last equality we used the definition $\check{Z}^{(b)}(\check{R}^{(b)};\hat{x}_0,\hat{x}_1,\ldots)=0$ of $\check{R}^{(b)}$.
	Similarly
	\begin{align*}
	D^{\mathrm{dil}} \check{R}^{(b)} &\stackrel{\text{Lem.~\ref{lem:Rder}}}{=} \left( I(b,1;\check{R}^{(b)}) - I(b,0;\check{R}^{(b)})  - \sum_{\ell=0}^{\infty} \hat{x}_{\ell}I(b,\ell;\check{R}^{(b)}) \right) \frac{\partial\check{R}^{(b)}}{\partial {\hat{x}_{b}}}\\
	&\mkern-8mu\stackrel{\text{Lem.~\ref{thm:IJprop}\ref{item:IJrel2}}}{=} \left( J(b;\check{R}^{(b)}) - \sum_{\ell=0}^{\infty} \hat{x}_{\ell}I(b,\ell;\check{R}^{(b)}) \right) \frac{\partial\check{R}^{(b)}}{\partial {\hat{x}_{b}}} = 0,
	\end{align*}
	which vanishes by the definition of $\check{R}^{(b)}$.
\end{proof}

\begin{lemma}\label{lem:strdilM}
For any $p,b\geq 0$ the power series $\check{M}_p^{(b)}(\hat{x}_0,\hat{x}_1,\ldots)$ satisfies
$D^{\mathrm{str}} \check{M}_p^{(b)} = 0$ and 
$D^{\mathrm{dil}} \check{M}_p^{(b)}= -\check{M}_p^{(b)}$.
\end{lemma}
\begin{proof}
First we rewrite the analogue of \eqref{eq:mhat} by introducing an additional formal variable $s$ via
\begin{align}
\check{M}_p^{(b)} &= Q_p(b,(1+r)\, \partial_r)\,(1+r)^{-b} \, \check{Z}^{(b)}(r; \hat{x}_{0},\hat{x}_{1}, \ldots) \Big|_{r=\check{R}^{(b)}}\nonumber\\
&= Q_p(b,(1+s)\, \partial_s)\, (1+s)^{-b}(1+\check{R}^{(b)})^{-b} \check{Z}^{(b)}(\,(1+s)(1+\check{R}^{(b)})-1; \hat{x}_{0},\hat{x}_{1}, \ldots) \Big|_{s=0}.\label{eq:Mcheckexpr}
\end{align}
Applying Lemma \ref{thm:IJprop} to the definition \eqref{eq:Zcheckdef} we easily find that
\begin{align}
D^{\mathrm{str}}\check{Z}^{(b)}(r;\hat{x}_0,\ldots) &= -(1+r)^{b+1} \frac{\partial}{\partial r}\left((1+r)^{-b}\check{Z}^{(b)}(r;\hat{x}_0,\ldots)\right), \label{eq:dstrzhat}\\
D^{\mathrm{dil}}\check{Z}^{(b)}(r;\hat{x}_0,\ldots) &= -\check{Z}^{(b)}(r;\hat{x}_0,\ldots).\label{eq:ddilzhat}
\end{align}
Therefore, with the help of Lemma \ref{lem:strdilR}, we have the identity
\begin{align*}
D^{\mathrm{str}}&\, (1+s)^{-b}(1+\check{R}^{(b)})^{-b} \check{Z}^{(b)}(\,(1+s)(1+\check{R}^{(b)})-1; \hat{x}_{0},\hat{x}_{1}, \ldots)\\
&=\left[(D^{\mathrm{str}}\check{R}^{(b)})\, (1+s)\frac{\partial}{\partial r}\left( (1+r)^{-b}\check{Z}^{(b)}(r;\hat{x}_0,\ldots)\right) +  (1+r)^{-b}D^{\mathrm{str}}\check{Z}^{(b)}(r;\hat{x}_0,\ldots)\right]_{r=(1+s)(1+\check{R}^{(b)})-1}\\
&\mkern-11mu\stackrel{\text{Lem.~\ref{lem:strdilR}}}{=}\left[(1+r)\frac{\partial}{\partial r}\left( (1+r)^{-b}\check{Z}^{(b)}(r;\hat{x}_0,\ldots)\right) +  (1+r)^{-b}D^{\mathrm{str}}\check{Z}^{(b)}(r;\hat{x}_0,\ldots)\right]_{r=(1+s)(1+\check{R}^{(b)})-1}
\end{align*}
of formal power series in $s,\hat{x}_0,\hat{x}_1,\ldots$. 
As a consequence of \eqref{eq:dstrzhat} it vanishes and therefore the same is true for $D^{\mathrm{str}}$ applied to $\check{M}_p^{(b)}$ given in \eqref{eq:Mcheckexpr}.
Similarly, since $D^{\mathrm{dil}}\check{R}^{(b)}=0$, we have
\begin{align*}
D^{\mathrm{dil}} &(1+s)^{-b}(1+\check{R}^{(b)})^{-b} \check{Z}^{(b)}(\,(1+s)(1+\check{R}^{(b)})-1; \hat{x}_{0},\hat{x}_{1}, \ldots) \\
& =  (1+r)^{-b}D^{\mathrm{dil}}\check{Z}^{(b)}(r;\hat{x}_0,\ldots)\Big|_{r=(1+s)(1+\check{R}^{(b)})-1} \\
&\stackrel{\eqref{eq:ddilzhat}}{=} -(1+s)^{-b}(1+\check{R}^{(b)})^{-b}\check{Z}^{(b)}(\,(1+s)(1+\check{R}^{(b)})-1;\hat{x}_0,\hat{x}_1,\ldots).
\end{align*}
This entails that $\check{M}_p^{(b)}$ from \eqref{eq:Mcheckexpr} satisfies $D^{\mathrm{dil}} \check{M}_p^{(b)}=-\check{M}_p^{(b)}$.
\end{proof}

\begin{proposition}\label{thm:stringdilatonhighergenus}
	The polynomials $\hat{N}^{(b)}_{g,n}$ satisfy for any $g,n\geq 1$ and $b,\ell_1,\ldots,\ell_n\geq 0$ the identities 
	\begin{align*}
	&\hat{N}^{(b)}_{g,n+1}(\ell_1,\ldots,\ell_n,1) = \sum_{j=1}^n \Big[-\ell_j\hat{N}^{(b)}_{g,n}(\ell_1,\ldots,\ell_n) + \sum_{k=1}^{\ell_j} 2 k\, \hat{N}^{(b)}_{g,n}(\ell_1,\ldots,\ell_{j-1},k,\ell_{j+1},\ldots,\ell_n) \\
	&\mkern355mu- \sum_{k=1}^{b} 2 k\, \hat{N}^{(b)}_{g,n}(\ell_1,\ldots,\ell_{j-1},k,\ell_{j+1},\ldots,\ell_n)\Big], \\
	&\hat{N}_{g,n+1}^{(b)}(\ell_1,\ldots,\ell_n,1) - \hat{N}_{g,n+1}^{(b)}(\ell_1,\ldots,\ell_n,0) = (n+2g-2)\, \hat{N}_{g,n}^{(b)}(\ell_1,\ldots,\ell_n).
	\end{align*}
\end{proposition}
\begin{proof}
Combining our definition of $\check{F}_g^{(b)}$ in terms of the moments $\check{M}_p^{(b)}$ with Lemma \ref{lem:strdilM}, we easily find
\begin{equation*}
D^{\mathrm{str}} \check{F}^{(b)}_1 = -\frac{1}{12}D^{\mathrm{str}} \log \check{M}_0 = 0, \qquad
D^{\mathrm{str}} \check{F}^{(b)}_g = D^{\mathrm{str}}P_g\left( \frac{1}{\check{M}_0}, \frac{\check{M}_1}{\check{M}_0},\ldots, \frac{\check{M}_{3g-3}}{\check{M}_0}\right) = 0\quad (g\geq 2).
\end{equation*}
Hence, for any $g,n\geq 1$ evaluating \eqref{eq:NfromFcheck} at $\ell_{n+1}=1$ gives
\begin{align*}
\hat{N}^{(b)}_{g,n+1}(\ell_1,\ldots,\ell_n,1) &= \frac{\partial^n}{\partial \hat{x}_{\ell_1}\cdots \partial \hat{x}_{\ell_n}} \frac{\partial}{\partial \hat{x}_1}\check{F}^{(b)}_{g}\Big|_{\hat{x}_0=\hat{x}_1=\ldots=0}\\
&= \frac{\partial^n}{\partial \hat{x}_{\ell_1}\cdots \partial \hat{x}_{\ell_n}}\sum_{\ell=0}^\infty \hat{x}_{\ell}\left( - \ell \frac{\partial}{\partial \hat{x}_{\ell}} +\sum_{k=1}^{\ell} 2k \frac{\partial}{\partial \hat{x}_{k}}-\sum_{k=1}^{b} 2k \frac{\partial}{\partial \hat{x}_{k}}\right)\hat{F}^{(b)}_{g}\Bigg|_{\hat{x}_0=\hat{x}_1=\ldots=0}\\
&= \sum_{j=1}^n \Big[-\ell_j\hat{N}^{(b)}_{g,n}(\ell_1,\ldots,\ell_n) + \sum_{k=1}^{\ell_j} 2 k\, \hat{N}^{(b)}_{g,n}(\ell_1,\ldots,\ell_{j-1},k,\ell_{j+1},\ldots,\ell_n) \\
&\mkern210mu- \sum_{k=1}^{b} 2 k\, \hat{N}^{(b)}_{g,n}(\ell_1,\ldots,\ell_{j-1},k,\ell_{j+1},\ldots,\ell_n)\Big].
\end{align*}
On the other hand $D^{\mathrm{dil}} \check{F}^{(b)}_1 = -\frac{1}{12}D^{\mathrm{dil}}\log \check{M}_0 = \frac{1}{12}$, while for $g\geq 2$ we find with the help of \eqref{eq:Ptilde} that
\begin{align*}
D^{\mathrm{dil}} \check{F}^{(b)}_g &= -D^{\mathrm{dil}} \check{M}_0^{2-2g}\tilde{P}_g\left(\frac{\check{M}_1}{\check{M}_0},\ldots, \frac{\check{M}_{3g-3}}{\check{M}_0}\right) \\
&=- (2g-2)\check{M}_0^{2-2g}\tilde{P}_g\left(\frac{\check{M}_1}{\check{M}_0},\ldots, \frac{\check{M}_{3g-3}}{\check{M}_0}\right) \\
&= (2g-2) \check{F}^{(b)}_g - (2g-2) \tilde{P}_g(0,\ldots,0).
\end{align*}
Finally, for any $g,n\geq 1$ we obtain
\begin{align*}
\hat{N}^{(b)}_{g,n+1}(\ell_1,\ldots,\ell_n,1)&-\hat{N}^{(b)}_{g,n+1}(\ell_1,\ldots,\ell_n,0)= \frac{\partial^n}{\partial \hat{x}_{\ell_1}\cdots \partial \hat{x}_{\ell_n}}\left( \frac{\partial}{\partial \hat{x}_1}-\frac{\partial}{\partial \hat{x}_0}\right)\check{F}^{(b)}_{g}\Big|_{\hat{x}_0=\hat{x}_1=\ldots=0}\\
&= \frac{\partial^n}{\partial \hat{x}_{\ell_1}\cdots \partial \hat{x}_{\ell_n}}\left( 2g-2 + \sum_{\ell=0}^\infty \hat{x}_{\ell}\frac{\partial}{\partial \hat{x}_{\ell}} \right)\check{F}^{(b)}_{g}\Big|_{\hat{x}_0=\hat{x}_1=\ldots=0}\\
&=(2g-2+n)\hat{N}^{(b)}_{g,n}(\ell_1,\ldots,\ell_n).
\end{align*}
\end{proof}

\subsection{Proof of Theorem \ref{thm:mainstringdilaton}}

The string and dilaton equations have been proven in Proposition \ref{thm:stringdilatonplanar} for $g=0$ and Proposition \ref{thm:stringdilatonhighergenus} for $g\geq 1$.
The polynomials $\hat{N}^{(b)}_{0,3}(\ell_1)$ and $\hat{N}^{(b)}_{g,1}(\ell_1)$ for $g\geq 1$ are independent of $b$ because of the remark just below Lemma \ref{lem:irrcriterion}: the criterion for being essentially $2b$-irreducible only puts restrictions on the lengths of essentially simple cycles enclosing at least two faces, of which there are none in $\hat{M}_{0,3}$ and $\hat{M}_{g,1}$.
A straightforward computation shows that $\hat{N}_{0,3}^{(b)}(\ell_1,\ell_2,\ell_3)=1$ and $\hat{N}_{1,1}^{(b)}(\ell_1)=(\ell_1^2-1)/12$.
The fact that the string and dilaton equations uniquely determine the symmetric polynomials $\hat{N}^{(b)}_{0,n}$ and $\hat{N}^{(b)}_{1,n}$ for $n>1$ in terms of $\hat{N}_{0,3}^{(b)}$ and $\hat{N}_{1,1}^{(b)}$ follows from the same reasoning as in \cite[Section 4.1]{Do2009}.\qed

\end{document}